\documentclass[11pt, a4paper]{amsart}
\usepackage{amsfonts,amssymb,amsmath, amsthm}
\usepackage{color}
\usepackage{mathrsfs}
\usepackage{mathtools}
\usepackage{lmodern}
\usepackage{latexsym}
\usepackage{bm}
\usepackage{enumerate}

\newcommand{\ddiv}{\operatorname{div}}

\allowdisplaybreaks

\theoremstyle{plain}
\newtheorem{theorem}{Theorem}[section]

\newtheorem{corollary}[theorem]{Corollary}

\newtheorem{proposition}[theorem]{Proposition}
\newtheorem{lemma}[theorem]{Lemma}

\theoremstyle{definition}
\newtheorem{definition}[theorem]{Definition}

\theoremstyle{remark}

\newtheorem{case}{Case}

\newtheorem*{claim*}{Claim}

\newtheorem{remark}[theorem]{Remark}

\numberwithin{equation}{section}

\def\eqn#1$$#2$${\begin{equation}\label#1#2\end{equation}}
\def\charfn_#1{{\raise1.2pt\hbox{$\chi_{\kern-1pt\lower3pt\hbox{{$\scriptstyle#1$}}}$}}}
\newcommand{\rif}[1]{(\ref{#1})}
\newcommand{\trif}[1] {\textnormal{\rif{#1}}}

\newcommand{\vp}{\varphi}

\newcommand{\e}{\varepsilon}

\DeclareMathOperator*{\osc}{osc}

\newcommand{\divo}{\operatorname{div}}

\newcommand{\supp}{\operatorname{supp}}

\def \R  {{\mathbb R}}

\def\er{\mathbb R}

\def\O{\mathcal{O}}

\def\YYint#1#2#3{{\setbox0=\hbox{$#1{#2#3}{\iint}$}
    \vcenter{\hbox{$#2#3$}}\kern-.58\wd0}}

 \newcommand{\means}[1]{-\hskip-1.00em\int_{#1}}

\def\mean#1{\mathchoice%
          {\mathop{\kern 0.2em\vrule width 0.6em height 0.69678ex depth -0.58065ex
                  \kern -0.8em \int}\limits_{\kern -0.4em#1}}%
          {\mathop{\kern 0.1em\vrule width 0.5em height 0.69678ex depth -0.60387ex
                  \kern -0.6em \int}\limits_{#1}}%
          {\mathop{\kern 0.1em\vrule width 0.5em height 0.69678ex
              depth -0.60387ex
                  \kern -0.6em \int}\limits_{#1}}%
          {\mathop{\kern 0.1em\vrule width 0.5em height 0.69678ex depth -0.60387ex
                  \kern -0.6em \int}\limits_{#1}}}

                  \def\mmean#1{\mathchoice%
          {\mathop{\kern 0.2em\vrule width 0.6em height 0.69678ex depth -0.58065ex
                  \kern -0.8em \iint}\limits_{\kern -0.4em#1}}%
          {\mathop{\kern 0.1em\vrule width 0.5em height 0.69678ex depth -0.60387ex
                  \kern -0.6em \iint}\limits_{#1}}%
          {\mathop{\kern 0.1em\vrule width 0.5em height 0.69678ex
              depth -0.60387ex
                  \kern -0.6em \iint}\limits_{#1}}%
          {\mathop{\kern 0.1em\vrule width 0.5em height 0.69678ex depth -0.60387ex
                  \kern -0.6em \iint}\limits_{#1}}}

\def\vintslides_#1{\mathchoice%
          {\mathop{\kern 0.1em\vrule width 0.5em height 0.697ex depth -0.581ex
                  \kern -0.6em \intop}\nolimits_{\kern -0.4em#1}}%
          {\mathop{\kern 0.1em\vrule width 0.3em height 0.697ex depth -0.604ex
                  \kern -0.4em \intop}\nolimits_{#1}}%
          {\mathop{\kern 0.1em\vrule width 0.3em height 0.697ex depth -0.604ex
                  \kern -0.4em \intop}\nolimits_{#1}}%
          {\mathop{\kern 0.1em\vrule width 0.3em height 0.697ex depth -0.604ex
                  \kern -0.4em \intop}\nolimits_{#1}}}

\newcommand{\aveint}[2]{\mathchoice%
          {\mathop{\kern 0.2em\vrule width 0.6em height 0.69678ex depth -0.58065ex
                  \kern -0.8em \intop}\nolimits_{\kern -0.45em#1}^{#2}}%
          {\mathop{\kern 0.1em\vrule width 0.5em height 0.69678ex depth -0.60387ex
                  \kern -0.6em \intop}\nolimits_{#1}^{#2}}%
          {\mathop{\kern 0.1em\vrule width 0.5em height 0.69678ex depth -0.60387ex
                  \kern -0.6em \intop}\nolimits_{#1}^{#2}}%
          {\mathop{\kern 0.1em\vrule width 0.5em height 0.69678ex depth -0.60387ex
                  \kern -0.6em \intop}\nolimits_{#1}^{#2}}}

                  \def\Xint#1{\mathchoice
{\XXint\displaystyle\textstyle{#1}}%
{\XXint\textstyle\scriptstyle{#1}}%
{\XXint\scriptstyle\scriptscriptstyle{#1}}%
{\XXint\scriptscriptstyle%
\scriptscriptstyle{#1}}%
\!\int}
\def\XXint#1#2#3{{\setbox0=\hbox{$#1{#2#3}{%
\int}$ }
\vcenter{\hbox{$#2#3$ }}\kern-.6\wd0}}
\def\barint{\,\Xint -} 
\def\bariint{\barint_{} \kern-.4em \barint}
\def\bariiint{\bariint_{} \kern-.4em \barint}
\renewcommand{\iint}{\int_{}\kern-.34em \int} 
\renewcommand{\iiint}{\iint_{}\kern-.34em \int} 

\newtoks\by
\newtoks\paper
\newtoks\book
\newtoks\jour
\newtoks\yr
\newtoks\pages
\newtoks\vol
\newtoks\publ
\def\et{ \& }
\def\name[#1, #2]{#1 #2}
\def\ota{{\hbox{\bf ???}}}
\def\cLear{\by=\ota\paper=\ota\book=\ota\jour=\ota\yr=\ota
\pages=\ota\vol=\ota\publ=\ota}
\def\endpaper{\the\by, \textit{\the\paper},
{\the\jour} \textbf{\the\vol} (\the\yr), \the\pages.\cLear}
\def\endbook{\the\by, \textit{\the\book},
\the\publ, \the\yr.\cLear}
\def\endpap{\the\by, \textit{\the\paper}, \the\jour.\cLear}
\def\endproc{\the\by, \textit{\the\paper}, \the\book, \the\publ,
\the\yr, \the\pages.\cLear}

\renewcommand{\d}{\, \mathrm{d}}

\title[Higher order interpolative geometries and evolutionary obstacle problems]{Higher order interpolative geometries and gradient regularity in evolutionary obstacle problems}

\author[Kim]{Sunghan Kim}
\address{Sunghan Kim\\ Department of Mathematics, Uppsala University, S-751 06 Uppsala, Sweden}
\email{sunghan.kim@math.uu.se}

\author[Nystr\"om]{Kaj Nystr\"om}
\address{Kaj Nystr\"om\\ Department of Mathematics, Uppsala University, S-751 06 Uppsala, Sweden}
\email{kaj.nystrom@math.uu.se}

\thanks{S.K. was supported by a grant from the Verge Foundation and K.N. was partially supported by grant 2022-03106 from the Swedish
research council (VR)}

\begin{document}

\maketitle

\date{today}

\begin{abstract}
\noindent We prove new optimal $C^{1,\alpha}$ regularity results for obstacle problems involving evolutionary $p$-Laplace type operators in the degenerate regime $p > 2$. Our main results include the optimal regularity improvement at free boundary points in  intrinsic backward $p$-paraboloids, up to the critical exponent, $\alpha \leq 2/(p-2)$,  and the optimal regularity across the free boundaries in the full cylinders up to a universal threshold. Moreover, we provide an intrinsic criterion by which the optimal regularity improvement at free boundaries can be extended to the entire cylinders. An important feature of our analysis is that we do not impose any assumption on the time derivative of the obstacle. Our results are formulated in function spaces associated to what we refer to as higher order or $C^{1,\alpha}$ intrinsic interpolative geometries.  \\

\noindent {\it Résumé}. Nous démontrons de nouveaux résultats optimaux de régularité $C^{1,\alpha}$ pour des problèmes d'obstacles impliquant des opérateurs de type $p$-Laplace évolutifs dans le régime dégénéré $p > 2$. Nos principaux résultats incluent l'amélioration optimale de la régularité aux points de frontière libre dans les $p$-paraboloides intrinsèques rétrogrades, jusqu'à l'exposant critique, $\alpha \leq 2/(p-2)$, et la régularité optimale à travers les frontières libres dans les cylindres complets jusqu'à un seuil universel. De plus, nous proposons un critère intrinsèque par lequel l'amélioration optimale de la régularité aux frontières libres peut être étendue à l'ensemble des cylindres. Une caractéristique importante de notre analyse est que nous n'imposons aucune hypothèse sur la dérivée temporelle de l'obstacle. Nos résultats sont formulés dans des espaces fonctionnels associés à ce que nous appelons des géométries interpolatives intrinsèques de premier ordre ou $C^{1,\alpha}$ plus élevées.\\


\end{abstract}




\section{Introduction}
This paper is devoted to the study of gradient regularity of solutions to obstacle problems involving quasilinear parabolic operators of the type
\eqn{generalH}
$$
  H u\equiv \divo  a(Du) - u_t,
$$
in cylindrical  domains of the form $
\O = \Omega \times (0,T) \subset \er^n \times \er
$, where $\Omega \subset \er^n$ is a bounded
Lipschitz domain, $T>0$ and $n \geq 2$.
The vector field $a \colon  \er^n \to \er^n$ is assumed to be $C^1$ regular and
satisfying the  {\em growth and ellipticity conditions}
\begin{equation}\label{asp}
\begin{dcases}
    |a(z)| +  (|z|^2 + s^2)^{\frac{1}{2}}  |\partial_z a(z)| \leq L(|z|^2+s^2)^{\frac{p-1}{2}},\\
   \langle \partial_z a(z)\xi, \xi\rangle \geq \nu(|z|^2+s^2)^{\frac{p-2}{2}}|\xi|^2,
\end{dcases}
\end{equation}
whenever $z, \xi \in \er^n$. Here $0< \nu\leq L$ and $s\geq 0$ are fixed parameters. Occasionally, we shall also require $a$ to be $C^2$ regular on $\mathbb R^n\setminus\{0\}$ and satisfying
\begin{equation}\label{eq:a-C2}
|\partial_z^2 a(z)| \leq L' (|z|^2 + s^2)^{\frac{p-3}{2}},\, z\in \mathbb R^n\setminus\{0\},
\end{equation}
for some fixed parameter $L' > 0$. Throughout the paper we will assume $p > 2$. The prototype for the operators considered is the $p$-parabolic, or evolutionary $p$-Laplace,  operator
\begin{equation}\label{modello}
\Delta_pu-u_t\equiv \ddiv (|Du|^{p-2} Du) - u_t.
\end{equation}

In the following we let  $
\partial_{\rm{P}}\O \equiv  (\bar \Omega \times \{0\}) \cup (\partial \Omega \times [0,T] )
$
denote the parabolic boundary of $\O$ and we will often refer to $\O$ as a parabolic domain.  Given a continuous boundary datum $g:\partial_{\rm P}\O\to\R$ and a continuous obstacle $\psi:\bar \O\to\R$ such that $g\geq\psi$ on $\partial_{\rm P} \O$, we  consider the obstacle problem
\begin{equation}\label{e-obs}
  \begin{dcases}
    \max\{H u, \psi-u\}=0 & \text{in}\ \O, \\
    u= g  & \text{on}\ \partial_{\rm P} \O.
  \end{dcases}
\end{equation}
We are particularly interested in the optimal regularity of the solution $u$ conditioned on the regularity of $g$ and $\psi$.  More precisely, the purpose of this paper is to establish {\it (optimal) regularity for $Du$ under stronger regularity assumptions on $\psi$ and $g$ beyond Lipschitz regularity in space}. In fact, to the best of our knowledge, {\it this is the first paper where optimal H\"older estimates for the gradient of solutions to $p$-parabolic obstacle problems are established}. The key novelty of our work arises in the introduction of the intrinsic geometry that reflects the (non-)degeneracy of the solution at reference points.

To formulate our results it is inevitable to introduce some notation and our intrinsic functional setting. It is a well established fact that due to the lack of homogeneity of the evolutionary $p$-Laplace operator and its generalizations,  the cylinders in which regularity is proved  have to depend intrinsically on the solution itself. The use of intrinsic geometry was pioneered in the work of E. DiBenedetto, see \cite{DBb}, and the philosophy is that for the regularity problem at hand, one has to unravel what cylinders to use, and one has to tailor the cylinders to the regularity one is aiming for. As a result, any precise regularity theory for quasilinear parabolic operators of $p$-Laplace type becomes inherently more complex compared to the corresponding (linear) theory in the case $p=2$.

\subsection{Modulus of continuity}
Throughout the paper, we will consider a  modulus of continuity $\omega: [0,\infty)\to [0,\infty)$, which is a concave and nondecreasing function satisfying $\omega(0) = 0$ and $\omega(1) = 1$. More precisely,  we will assume that there are constants $A\geq 1$ and {$0 < \alpha \leq\frac{2}{p-2}$} such that
\begin{equation}\label{eq:omega prop}
1\leq \frac{\omega(\varrho)}{\omega(r)}\leq A\left(\frac{\varrho}{r}\right)^\alpha,\quad \text{whenever}\quad \varrho \geq r.
\end{equation}
The growth condition in \eqref{eq:omega prop} obviously allows for H\"older moduli of continuity, i.e., $\omega(r) = r^\alpha$, but the condition also allows for logarithmic type moduli of continuity which do not even satisfy the Dini condition, e.g., $\omega(r) = \max\{ 1/\log(e/r) , 1\}$.

\subsection{Intrinsic space-time cylinders and $p$-paraboloids}
Let $B_r$ denote the ball in $\R^n$ with center at the origin and with radius $r$. Given $r$ and $\mu\geq 0$, we introduce the backward and forward intrinsic space-time cylinders centered at (0,0),
\begin{equation}\label{eq:Cyl-}
\begin{aligned}
Q_r^-(\mu) \equiv  B_r \times (-\mu^{2-p}r^2, 0), \quad Q_r^+(\mu) \equiv  B_r \times (0,\mu^{2-p}r^2),
\end{aligned}
\end{equation}
as well as the entire cylinders
$$
Q_r (\mu) \equiv   B_r \times (-\mu^{2-p}r^2 , \mu^{2-p}r^2).
$$
For brevity, we shall write $Q_r^-$, $Q_r^+$, and $Q_r$, for  the cylinders $Q_r^-(1)$, $Q_r^+(1)$,  and  $Q_r(1)$, respectively. We also introduce the backward $p$-paraboloids
\begin{equation}\label{eq:p-para}
P_r^-(\mu) \equiv  \left\{ (x,t) \in Q_r^-(\mu) : \frac{|x|^p}{-t} \leq \mu^{p-2}\right\}.
\end{equation}
Regularity will be measured in intrinsic cylinders and $p$-paraboloids of the type
\begin{align}\label{cyla}
(x_0,t_0) + Q_r^\pm ( \lambda\vp(r) ),
\end{align}
and
\begin{equation}\label{eq:para}
(x_0,t_0) + P_r^- (\vp(r)),
\end{equation}
respectively,  where $\vp(r)$ denotes a carefully chosen intrinsic scaling factor (in time).

\subsection{Higher order intrinsic interpolative geometries}  In Definitions \ref{def:point-reg}, \ref{def:point-regaa} and \ref{def:full-reg} below we introduce the
$C^{1,\alpha}$ intrinsic function spaces $\tilde C_{s}^{1,\omega}$ and $\tilde C_{s,{\rm par}}^{1,\omega}$. These spaces of functions are defined in terms of pointwise approximations by time-independent affine functions. Our key estimates will be built on these spaces. The definitions, which incarnate the  higher order intrinsic interpolative geometries, are quite involved due the intrinsic geometry, so it may at first glance be difficult to fully grasp their meaning. For this reason, we here give a brief informal discussion of the higher order intrinsic interpolative geometries which lead up to the definitions.

Let $u$ be a continuous function on a parabolic domain $\O$ in $\R^{n+1}$, let $\omega$ be a modulus of continuity, and let $s$ be a nonnegative constant. Assume, for the sake of simplicity, that $$(x_0,t_0) + Q_{1}\subset \O,$$ and consider all time-independent affine functions $\ell$ such that
\begin{equation}\label{Introoo1}
\mbox{$\ell(x_0) = u(x_0,t_0)$ and  $|D\ell| + s \leq 1$}.
\end{equation} In this setting, our intrinsic function spaces are stated, at  $(x_0,t_0)\in\O$ and for $r\leq 1$, in terms of the existence of a time-independent affine functions $\ell$ satisfying \eqref{Introoo1} and  such that
$|u - \ell|/(r\omega(r))$ is bounded by $\lambda$, measured over the cylinder $(x_0,t_0) + Q_r(\lambda\vp(r))$, where
\begin{align*}\vp(r) = \vp(r;x_0,t_0) \equiv  \max\{\omega(r), |D\ell| + s \},\ r>0.
\end{align*}
A posteriori, $Du(x_0,t_0)$  exists and  $D\ell=Du(x_0,t_0)$. Let us here assume $\omega(r)=r^\alpha$, {$\alpha\leq\frac{2}{p-2}$}.

We start by considering two extremal cases. First, if $|D\ell| + s=0$, then  $\vp(r) = \omega(r)$ for $r\in (0,1)$ and
\begin{align}\label{exp2}(x_0,t_0) + Q_r(\lambda\vp(r))=(x_0,t_0)+B_r \times (-\lambda^{2-p}r^{2-(p-2)\alpha} , \lambda^{2-p}r^{2-(p-2)\alpha}).
\end{align}
In this case, our function space measures the deviation of $u$ from $u(x_0,t_0) $ according to \begin{equation}
    \sup_{(x_0,t_0) + Q_r(\lambda \omega(r))}
    |u - u(x_0,t_0) |
    \leq
    \lambda r^{1+\alpha},\  \forall r\in (0,1).
\end{equation}
Second, if $|D\ell| + s=1$, then $\vp(r) = 1$ for $r\in (0,1)$ and
\begin{align}\label{exp3}(x_0,t_0) + Q_r(\lambda\vp(r))=(x_0,t_0)+B_r \times (-\lambda^{2-p}r^2 , \lambda^{2-p}r^2).
\end{align}
The cylinders in \eqref{exp3} are the cylinders that naturally appear in the study of uniformly elliptic second order parabolic equations. In this case, our function space measures the deviation of $u$ from $\ell$ according to \begin{equation}
    \sup_{(x_0,t_0) + Q_r(\lambda)}
    |u - \ell |
    \leq
    \lambda r^{1+\alpha},\ \forall r\in (0,1).
\end{equation}
In both these extremal cases $Du(x_0,t_0)$  exists,  $D\ell=Du(x_0,t_0)$, and  $\lambda$ is a measurement of $C^{1,\alpha}$.

More generally, and going beyond the extremal cases,  the function spaces we introduce interpolate between the degenerate case, $|D\ell| + s\approx 0$, and the nondegenerate case, $|D\ell| + s\approx 1$. Indeed, given a threshold $\mu\equiv |D\ell| + s\in (0,1)$, we introduce the (smallest degenerate) scale $r_\omega\equiv \omega^{-1}(\mu)=\mu^{1/\alpha}$, and we let
 \begin{align}\label{exp4}\vp(r) = \vp(r;x_0,t_0) = \chi_{[r_\omega,1)}(r)\omega(r)+\chi_{(0,r_\omega)}(r)\mu=\chi_{[\mu^{1/\alpha},1)}(r)r^\alpha+\chi_{(0,\mu^{1/\alpha})}(r)\mu,\ r>0,
\end{align}
where $\chi_I$ is the indicator function for an interval $I\subset(0,1)$. Given $\vp$ as in \eqref{exp4} we then consider the cylinders $(x_0,t_0) + Q_r(\lambda\vp(r))$ and we  measure the deviation of $u$ from $\ell$ according to \begin{equation}
    \sup_{(x_0,t_0) + Q_r(\lambda \vp(r))}
    |u - \ell |
    \leq
    \lambda r^{1+\alpha}.
\end{equation}
Hence, $\lambda$ is a measurement of $C^{1,\alpha}$ which takes the level of degeneracy of $|D\ell| + s=|D u| + s$ at $(x_0,t_0)$ into account by interpolating between the degenerate and the nondegenerate scales: this is the essence of our higher order intrinsic interpolative geometries underlying our function space $\tilde C_{s}^{1,\omega}$.

As it turns out, at free boundary points we will also consider measurements in $(x_0,t_0) + P_r^-(\vp(r))$, where again in the case of $\omega(r) \equiv r^\alpha$ we have two extreme cases
\begin{equation}\label{exp5}
P_r^-(\vp(r)) =
\begin{cases}
    \{ (x,t) \in B_r\times (-r^{2-(p-2)\alpha},0) : {|x|^p} \leq -tr^{\alpha(p-2)}\}, & \text{if }|D\ell| + s = 0, \\
    \{ (x,t) \in B_r\times (-r^2,0) : {|x|^p} \leq -t\}, & \text{if }|D\ell| + s = 1,
\end{cases}
\end{equation}
for every $r\in(0,1)$. In general, given a threshold  $\mu\equiv |D\ell| + s\in (0,1)$, $P_r^-(\vp(r))$ interpolates between the above two geometries. This intrinsic geometry underlies our function space $\tilde C_{s,{\rm par}}^{1,\omega}$ where the subscript par refer to paraboloid.

\subsection{Function spaces} We here give the formal definitions of the $C^{1,\alpha}$ intrinsic function spaces $\tilde C_{s}^{1,\omega}$ and $\tilde C_{s,{\rm par}}^{1,\omega}$ which incarnate the higher order intrinsic interpolative geometries.

\begin{definition}[Pointwise regularity]\label{def:point-reg}
Let $u$ be a continuous function on a parabolic domain $\O$ in $\R^{n+1}$, let $\omega$ be a modulus of continuity {with $\alpha \leq \frac{2}{p-2}$}, ${\bar r}>0$, $\bar\mu>0$, and let $s$ be a nonnegative constant. Assume that $$\mbox{$(x_0,t_0) + Q_{\bar r}(\bar\mu) \subset \O$,}$$
and consider all time-independent affine functions $\ell$ such that
\begin{equation}\label{Introoo2}
\mbox{$\ell(x_0) = u(x_0,t_0)$ and  $|D\ell| + s \leq \bar\mu$}.
\end{equation}
 Given $\ell$ we let
\begin{equation}\label{ttscale}
 \vp (r) \equiv  \vp (r\, ; x_0,t_0,{\bar r},\bar\mu) \equiv  \max\left\{ \omega(r), \frac{\omega({\bar r})}{\bar\mu}(|D\ell| + s) \right\},\ r>0.
 \end{equation}
We say that $u \in \tilde C_{s,\pm}^{1,\omega}(x_0,t_0;\bar r, \bar\mu)$, if there exists a finite constant $\lambda$ and a time-independent affine function $\ell$ as in \eqref{Introoo2}, such that
$$ \sup_{(x_0,t_0) + (Q_r^\pm(\lambda\vp(r))\cap Q_{\bar r}(\bar\mu))} |u - \ell | \leq \lambda r\omega(r),\,\forall r\in(0,{\bar r}),$$
with $\vp$ is as in \eqref{ttscale}. We let $[u]_{\tilde C_{s,\pm}^{1,\omega}(x_0,t_0; {\bar r},\bar\mu)}$ denote the smallest such $\lambda$. We say that $u\in \tilde C_s^{1,\omega}(x_0,t_0;{\bar r},\bar\mu)$ if $u\in \tilde C_{s,+}^{1,\omega}(x_0,t_0;{\bar r},\bar\mu)\cap \tilde C_{s,-}^{1,\omega}(x_0,t_0;{\bar r},\bar\mu)$ and we let
$$
[u]_{\tilde C_s^{1,\omega}(x_0,t_0;{\bar r},\bar\mu)} \equiv  \min\left\{[u]_{\tilde C_{s,+}^{1,\omega}(x_0,t_0;{\bar r},\bar\mu)} ,[u]_{\tilde C_{s,-}^{1,\omega}(x_0,t_0;{\bar r},\bar\mu)} \right\}.
$$
\end{definition}

To study estimates in backward $p$-paraboloids we introduce the following class of functions.
\begin{definition}[Pointwise regularity in backward $p$-paraboloids]\label{def:point-regaa}
Let $u$ be a continuous function on a parabolic domain $\O$ in $\R^{n+1}$, let $\omega$ be a modulus of continuity {with $\alpha\leq\frac{2}{p-2}$}, ${\bar r}>0$, $\bar\mu>0$, and let $s$ be a nonnegative constant. Assume that $$\mbox{$(x_0,t_0) + Q_{\bar r}(\bar\mu) \subset \O$.}$$
 We say that $u \in \tilde C_{s,{\rm par}}^{1,\omega}(x_0,t_0;{\bar r},\bar\mu)$  if there exists a finite constant $\lambda$ and a time-independent affine function $\ell$ as in \eqref{Introoo2}, such that
 $$
\sup_{(x_0,t_0) + (P_r^-(\omega(r))\cap Q_{\bar r}(\bar\mu))} |u - \ell| \leq \lambda r\omega(r),\,\forall r \in (r_\omega,{\bar r}),
$$
where $r_\omega := \inf\{ r \in (0,{\bar r} ] :\vp(r) = \omega(r)\} \in [0,{\bar r}]$ and where $\vp$ is defined as in \eqref{ttscale}. We let $[u]_{\tilde C_{s,{\rm par}}^{1,\omega}(x_0,t_0; {\bar r},\bar\mu)}$ denote the smallest such $\lambda$.
\end{definition}

Note that in Definition \ref{def:point-regaa}, the definition of $u\in \tilde C_{s,{\rm par}}^{1,\omega}(x_0,t_0;\bar r,\bar\mu)$ is stated in terms of measurements of  $u-\ell$ over the sets $(x_0,t_0) + P_r^-(\omega(r))$, and especially if $\omega(r) \equiv r^\alpha$, then
$$
P_r^-(\omega(r)) = \{ (x,t) \in B_r\times (-r^{2-(p-2)\alpha},0) : {|x|^p} \leq -tr^{\alpha(p-2)}\}.
$$
 When $\bar r = \bar \mu = 1$, then $r_\omega \equiv \inf\{ r \in (0,1]: \vp(r)=\omega(r)\} = \omega^{-1}(\mu)$, with $\mu = |D\ell| + s \in [0,1]$. We refer to $r_\omega$ as  the smallest degenerate scale. Note that for any degenerate scale $r \in (r_\omega,1)$ (recall $\bar r= \bar \mu = 1$), we have
\begin{eqnarray}\label{Replaceosc}
\left| \sup_{(x_0,t_0) + P_r^-(\omega(r))} |u -\ell| - \sup_{(x_0,t_0) + P_r^-(\omega(r))} |u-u(x_0,t_0)| \right| \leq r|D\ell| \leq r\omega(r).
\end{eqnarray}
We remark that the inequality in the last display may not be true for $r < r_\omega$ when $|D\ell| > 0$, and this is the reason we call $r_\omega$ the {\it smallest degenerate scale}. In particular,  in the definition of $\tilde C_{s,{\rm par}}^{1,\omega}(x_0,t_0;{\bar r},\bar\mu)$, the supremum of $|u-\ell|$ is measured only for the degenerate scales, and \eqref{Replaceosc} shows that we get a definition equivalent to Definition \ref{def:point-regaa}, up to a multiplicative constant, if we replace the supremum by $$\osc_{(x_0,t_0) + P_r^-(\omega(r))} u.$$
Here, and in general, given $A\subset\mathbb R^{n+1}$ and
a function $f:A\to \mathbb R^{m}$, $m\geq 1$, we let
$$\osc_{A} f:=\sup_{(x,t),(\tilde x,\tilde t)\in A}|f(x, t)-f(\tilde x,\tilde t)|,$$
denote the oscillation of $f$ on $A$.

To study estimates in the entire domain we introduce the following class of functions.

\begin{definition}[Full regularity]\label{def:full-reg}
Let $u$ be a continuous function on a parabolic domain $\O$ in $\R^{n+1}$, let $\omega$ be a modulus of continuity {with $\alpha\leq\frac{2}{p-2}$}, and let $s$ be a nonnegative constant. We say that $u \in \tilde C_s^{1,\omega}(\O)$, if there is a finite constant $\lambda$ such that if $(x_0,t_0)\in\O$, then there is a time-independent affine function $\ell$ such that
\begin{align}\label{equal1}
\sup_{((x_0,t_0) + Q_r(\lambda\vp(r)))\cap\O} |u - \ell| \leq \lambda r\omega(r),
\end{align}
for all $r > 0$ where $$\vp(r) = \vp(r;x_0,t_0) \equiv  \max\{\omega(r), |D\ell| + s \}.$$ We let $[u]_{\tilde C_s^{1,\omega}(\O)}$ denote the smallest $\lambda$ for which \eqref{equal1} holds uniformly for all $(x_0,t_0) \in \O$.
\end{definition}

We shall also use $C^{0,\omega}(\O)$ to denote the usual function space consisting of all continuous functions on $\O$ with modulus of continuity $\omega$, and $C^{k,\omega}(\O)$ to denote the usual function space consisting of all $k$-times continuously differentiable functions whose $k$-th order derivatives (in both space and time) belong to $C^{0,\omega}(\O)$.


\subsection{Main results}

We here collect our main results. Our first result concerns  interior H\"older estimates for the gradient of solutions to the obstacle problem in terms of the intrinsic geometry.

\begin{theorem}[Regularity in the interior]\label{thm:C1a-int}
Let $H$ be as in \eqref{generalH}, with $(a,s)$ as in \eqref{asp} and \eqref{eq:a-C2}, let $\psi$ be an obstacle in $\O$, and let $u$ be a weak solution to $\max\{H u, \psi-u\} = 0$ in $\O$. Then there exists $\alpha_h \equiv  \alpha_h(n,p,\nu,L)$, $\alpha_h \in(0,1)$, such that if $\omega$ verifies \eqref{eq:omega prop} with $\alpha \in (0,\alpha_h)$, and if
$$
|Du(x_0,t_0)| + s + \omega({\bar r}) [\psi]_{\tilde C_s^{1,\omega}(x_0,t_0;{\bar r},\bar\mu)} +  \frac{1}{{\bar r}}  |u - u(x_0,t_0)| \leq \bar\mu\quad\text{in }(x_0,t_0) +Q_{\bar r}(\bar\mu) \subset \O,
$$
for some ${\bar r}>0$ and $\bar\mu>1$, then $u\in \tilde C_s^{1,\omega}(x_0,t_0;{\bar r},\bar\mu)$, and
$$
[u ]_{\tilde C_s^{1,\omega}(x_0,t_0;{\bar r},\bar\mu)} \leq \frac{\lambda\bar\mu}{\omega({\bar r})},
$$
for some constant $\lambda \equiv  \lambda(n,p,\nu,L,L',A,\alpha)$.
\end{theorem}

{
To help the readers understand the intrinsic $C^{1,\omega}$-estimate above, let us restate the result in terms of the usual function space.
\begin{corollary}
Under the setting of Theorem \ref{thm:C1a-int}, there is a time-independent affine function $\ell$, satisfying $|D\ell| + s \leq \bar\mu$, such that for any $r\in (0,\bar r)$,
$$
\begin{dcases}
\sup_{(x_0,t_0) + Q_r(\bar\lambda\omega(r))} | u - \ell | \leq
\bar\lambda r\omega(r), & \text{if } \omega(r) \geq \frac{\omega(\bar r)}{\bar\mu}(|D\ell| + s), \\
\sup_{(x_0,t_0) + Q_r(\frac{\bar\lambda}{\bar\mu}\omega(\bar r))} | u - \ell | \leq \bar\lambda r\omega(r), & \text{if  } \omega(r) \leq \frac{\omega(\bar r)}{\bar\mu}(|D\ell| + s),
\end{dcases}
$$
with $\bar \lambda := \lambda\bar\mu(\omega(\bar r))^{-1}$. 
\end{corollary}
}

\begin{remark}\label{rem:C1a-int}
Note that the reference point $(x_0,t_0)$ in Theorem \ref{thm:C1a-int} does {\it not} necessarily lie in the coincidence set $\{ u = \psi\}$. In fact, if either $(x_0,t_0) + Q_{\bar r}(\bar\mu) \subset \{ u > \psi\}$ or $(x_0,t_0) \in \{ u = \psi\}$, then the assertion of Theorem \ref{thm:C1a-int} holds assuming only the basic structure condition \eqref{asp}, see Proposition \ref{prop:C1a-int} and Proposition \ref{prop:C1a}. The additional condition \eqref{eq:a-C2} is only imposed to allow us to combine two separate estimates in the intermediate case when  $(x_0,t_0) \in \{ u >\psi\}$, but $((x_0,t_0) + Q_{\bar r}(\bar\mu)) \cap\{u = \psi\} \neq\emptyset$. However, if we accept to lose some order of regularity of $u$ relative to that of $\psi$, then we are free to remove the additional assumption \eqref{eq:a-C2}, see Remark \ref{rem:C1a-gen} for more discussions. Finally, we remark that the prototype for our operators, the evolutionary $p$-Laplace operator, $\Delta_p u - u_t$, verifies both \eqref{asp} and \eqref{eq:a-C2}.
\end{remark}

\begin{remark} Note that in Theorem \ref{thm:C1a-int} we prove optimal regularity in the entire intrinsic cylinders, given the regularity of the obstacle,
up to the threshold $\alpha_h$ (the implicit interior regularity exponent), i.e., for less regular obstacles. As $\alpha_h$ is implicit, given the regularity of the obstacle,  the {\it optimal regularity in the full intrinsic cylinders} remains unknown; more precisely, the sharp regularity threshold up to which the estimate holds in the full intrinsic cylinders remains unknown. Theorem \ref{thm:C1a-int} (see also Proposition \ref{prop:C1a}) shows that the threshold has a lower bound $\alpha_h$, and the example with the traveling wave solutions outlined in Subsection \ref{sec:novel} below yields the upper bound ${1}/{(p-2)}$. It is an intriguing problem to capture the optimal threshold.
\end{remark}

\begin{remark} A consequence of the main result in \cite{KMNobs}, see  Corollary 1.1 in \cite{KMNobs}, is an endpoint result which states that if $D\psi$, $g\in \mathrm{L}^\infty(\O)$, then  any weak solution $u$ to \eqref{e-obs} satisfies $Du\in\mathrm{L}_{\text{loc}}^\infty(\O)$.
\end{remark}
An important feature in the obstacle problem for the (stationary) $p$-Laplace operator, the $p$-obstacle problem for short,  is the improvement of the regularity at free boundary points. That is, for the $p$-obstacle problem, the solution is as smooth as the obstacle, up to $C^{1,1}$, at every free boundary point, though the solution is,  independent of the smoothness of the obstacle, only of class $C^{1,\alpha}$, for some (small) universal $\alpha$, inside the continuation domain (set). In our evolutionary setting, we obtain the optimal improvement of regularity at free boundary points, yet only in the intrinsic backward $p$-paraboloids. In fact, the involvement of the $p$-paraboloids is sharp for regular obstacles, and we shall come back to this important matter later.  More specifically, we prove the following theorem.

\begin{theorem}[Regularity improvement at free boundaries]\label{thm:C1a}
Let $H$ be as in \eqref{generalH}, with $(a,s)$ as in \eqref{asp}, let $\psi$ be an obstacle in $\O$, and let $u$ be a solution to $\max\{ Hu, \psi-u \} = 0$ in $\O$. Suppose that $u(x_0,t_0) = \psi(x_0,t_0)$ at some $(x_0,t_0)\in \O$, that $\omega$ verifies \eqref{eq:omega prop} for some $\alpha \in (0,{2}/{(p-2)}]$, and that
$$
|D\psi(x_0,t_0)| + s  + \omega({\bar r}) [\psi]_{\tilde C_s^{1,\omega}(x_0,t_0;{\bar r},\bar\mu)} \leq \bar\mu \quad\text{and}\quad (x_0,t_0) +Q_{\bar r}(\bar\mu) \subset \O,
$$
for some ${\bar r} > 0$ and $\bar\mu \geq 1$. Then $u\in \tilde C_{s,{\rm par}}^{1,\omega}(x_0,t_0;{\bar r},\bar\mu)$ and
$$
[u ]_{\tilde C_{s,{\rm par}}^{1,\omega}(x_0,t_0;{\bar r},\bar\mu)} \leq \frac{\lambda\bar\mu}{\omega({\bar r})},
$$
for some $\lambda \equiv \lambda(n,p,\nu,L,A,\alpha)$.
\end{theorem}

{

\begin{corollary}
Under the setting of Theorem \ref{thm:C1a}, there is a time-independent affine function $\ell$, satisfying $|D\ell| + s \leq \bar\mu$, such that for any $r\in (0,\bar r)$,
$$
\sup_{(x_0,t_0) + P_r^-(\omega(r))} | u - \ell | \leq
 \bar\lambda r\omega(r), \quad\text{whenever } \omega(r) \geq \frac{\omega(\bar r)}{\bar\mu} (|D\ell| + s),
$$
with $\bar\lambda:= \lambda\bar\mu(\omega(\bar r))^{-1}$. 
\end{corollary}
}

\begin{remark}\label{rem:C1a} Note that we do {\it not} assume any upper bound of $u - \psi(x_0,t_0)$ in $(x_0,t_0) + Q_{\bar r}(\bar\mu)$ in the statement of
Theorem \ref{thm:C1a} (a lower bound is given by the regularity of the obstacle $\psi$). Thus, the regularity improvement is purely due to the presence of the obstacle. The restriction in the backward $p$-paraboloids is natural, and cannot be improved in general, and as mentioned we shall discuss this in more detail later, see Subsection \ref{sec:novel}. Roughly speaking, due to the lack of the strong minimum principle in our setting, one cannot expect universal estimates over entire cylinders that are independent of the boundary values.
\end{remark}
\begin{remark}\label{rem:C1a+} In the proof of Theorem \ref{thm:C1a} the time-independent affine function $\ell(x)\equiv \psi(x_0,t_0)+ D\psi(x_0,t_0)\cdot(x-x_0)$ is central.  Using this affine function, we note that  in the nondegenerate scenario $|D\ell| + s=|D\psi(x_0,t_0)| + s > 0$, the corresponding estimate in Theorem \ref{thm:C1a} stops at the scale $r_\omega := \inf\{ r \in (0,{\bar r} ] :\vp(r) = \omega(r)\} >0$ where the problem is no longer degenerate, see Definition \ref{def:point-regaa}. However, by incorporating the boundary values of $u - \psi(x_0,t_0)$ at the scale $r_\omega$ of the nondegeneracy\footnote{To be explicit, $r_\omega$ is the largest $r \in (0,{\bar r}]$ such that $\omega(r) \bar\mu \leq \omega({\bar r})(|D\psi(x_0,t_0)| + s)$.}, we manage to extend the optimal intrinsic estimate to the entire lower cylinders, see Lemma \ref{lem:C1a-ndeg}. This being said, we could have given an alternative definition of the space $\tilde C_{s,{\rm par}}^{1,\omega}(x_0,t_0;{\bar r},\bar\mu)$, and hence an alternative formulation of Theorem \ref{thm:C1a},  which captures the measurement of the supremum of $|u-\ell|$, over the backward paraboloids, for all scales $r\in (0,{\bar r})$ and not only for the degenerate scales $(r_\omega,{\bar r})$. However, a defect in such a formulation  is the loss of the uniform bound on the corresponding seminorm.
\end{remark}

 Theorem \ref{thm:C1a} shows that at free boundary points,  the regularity threshold can be raised to ${2}/{(p-2)}$ in the backward $p$-paraboloids $(P_{\bar r}^-)$. To extend this result to the full intrinsic cylinders it is reasonable to believe that the optimal regularity threshold may appear in an intrinsic way. In fact, the following theorem states that there is another intrinsic quantity such that if the quantity stays bounded, then the optimal regularity improvement holds in the intrinsic cylinders.

\begin{theorem}\label{thm:C1a-re} Let $H$ be as in \eqref{generalH}, with $(a,s)$ as in \eqref{asp}, let $\psi$ be an obstacle in $\O$, and let $u$ be a solution to $\max\{ Hu, \psi-u \} = 0$ in $\O$. Suppose that $u(x_0,t_0) = \psi(x_0,t_0)$ at some $(x_0,t_0)\in \O$, that $\omega$ verifies \eqref{eq:omega prop} for some $\alpha \in (0,{2}/{(p-2)}]$, and that
$$
|D\psi(x_0,t_0)| + s  + \omega({\bar r}) [\psi]_{\tilde C_s^{1,\omega}(x_0,t_0;{\bar r},\bar\mu)} + \frac{1}{{\bar r}} |u - u(x_0,t_0)| \leq \bar\mu\quad\text{in }(x_0,t_0) +Q_{\bar r}(\bar\mu) \subset \O,
$$
for some ${\bar r} > 0$ and $\bar\mu \geq 1$. If there exists a constant $\theta > 1$ such that
\begin{equation}\label{eq:C1a-re-asmp}
\frac{\inf\{ \lambda > 1: |u- u(x_0,t_0)|\leq \lambda r\omega(r)\text{ in }(x_0,t_0) + Q_r(\omega(r))\}}{\inf\{ \lambda > 1: |u- u(x_0,t_0)|\leq \lambda r\omega(r)\text{ in }(x_0,t_0)+Q_r(\lambda\omega(r))\}}\leq \theta,
\end{equation}
for every $r$ such that  $\omega({\bar r})(|D\psi(x_0,t_0)| + s)/\bar\mu<\omega(r)$, then
$$
[u ]_{\tilde C_s^{1,\omega}(x_0,t_0; {\bar r},\bar\mu)} \leq \frac{\lambda\bar\mu}{\omega({\bar r})},
$$
for some $\lambda\equiv\lambda(n,p,\nu,L,A,\alpha,\theta)$.
\end{theorem}

\begin{remark}\label{rem:C1a-re}
We emphasize that the criterion in \eqref{eq:C1a-re-asmp} is {\it not} redundant  for the conclusion concerning the optimal regularity in the intrinsic cylinders at free boundary points. This will be shown with an example in Subsection ~\ref{sec:novel} below. In fact, one can prove that the criterion in \eqref{eq:C1a-re-asmp} is both sufficient and necessary for the conclusion of Theorem \ref{thm:C1a-re}.
\end{remark}

\begin{remark}\label{rem:C1a-re2}
The criterion \eqref{eq:C1a-re-asmp} always holds (for some $\theta$) when $|D\psi(x_0,t_0)| + s > 0$, i.e., in the nondegenerate scenario. This justifies our comment in Remark \ref{rem:C1a+} that our solution belongs to the class $\tilde C_{s,{\rm par}}^{1,\omega}$ even if we extend the definition to the nondegenerate scales; but again, we lack a uniform bound on the corresponding norm, as we do not have a uniform control over $\theta$.
\end{remark}

At the initial (time) layer/state we obtain the optimal regularity in the forward intrinsic cylinders. Moreover, the assertion holds even when the reference point is {\it not} necessarily a contact point. Our proof relies on  uniqueness results for homogeneous Cauchy problems. In this context we prove the following theorem.

\begin{theorem}[Regularity at the initial layer]\label{thm:C1a-ini}
Let $H$ be as in \eqref{generalH}, with $(a,s)$ as in \eqref{asp}, let $\psi$ be an obstacle in $(x_0,t_0) + \overline{Q_{\bar r}^+(\bar\mu)}$, let $g$ be an initial datum with $g \geq \psi(\cdot,t_0)$ on $B_{\bar r}(x_0)$, and let $u$ be a continuous solution to
$$
  \begin{cases}
    \max\{H u, \psi-u\}=0 & \text{in}\ (x_0,t_0)+ Q_{\bar r}^+(\bar\mu), \\
    u= g  & \text{on}\ (x_0,t_0) + B_{\bar r}\times\{0\}.
  \end{cases}
$$
Suppose that $\omega$ verifies \eqref{eq:omega prop} with $\alpha \in (0,{2}/{(p-2)}]$, that
$$
|Dg(x_0)| + |D\psi(x_0,t_0)| + s + \frac{1}{{\bar r}}|u - g(x_0)| \leq \bar\mu\quad\text{in }(x_0,t_0) + \overline{Q_{\bar r}^+(\bar\mu)},
$$
for some ${\bar r} > 0$ and $\bar\mu \geq 1$, and that
$$
[g]_{C^{1,\omega}(x_0;{\bar r})} + [\psi]_{\tilde C_{s,+}^{1,\omega}(x_0,t_0;{\bar r},\bar\mu)} \leq \frac{\bar\mu}{\omega({\bar r})}.
$$
Then $u\in \tilde C_{s,+}^{1,\omega}(x_0,t_0;{\bar r},\bar\mu)$ and
$$
[u ]_{\tilde C_{s,+}^{1,\omega}(x_0,t_0;{\bar r},\bar\mu)} \leq \frac{\lambda\bar\mu}{\omega({\bar r})},
$$
for some $\lambda \equiv \lambda(n,p,\nu,L,A,\alpha)$.
\end{theorem}

One may also need a regularity theorem in the entire domain, rather than just the pointwise approximation results. That is the content of our last theorem.

\begin{theorem}[Regularity across free boundaries] \label{thm:C1a-gen}
Let $H$ be as in \eqref{generalH}, with $(a,s)$ as in \eqref{asp} and \eqref{eq:a-C2}, let $\omega$ verify \eqref{eq:omega prop} with $\alpha \in (0,\alpha_h)$, and let $u$ be a solution to \eqref{e-obs}. Then,  if $\psi \in \tilde C_s^{1,\omega}(\O)$, then $u \in \tilde C_s^{1,\omega}(\O')$ for every compact subset $\O'$ of $\O$. Moreover, if $g(\cdot,0) \in C^{1,\omega}(\Omega)$,  then $u \in \tilde C_s^{1,\omega}(\O')$ for every $\O'\equiv\Omega'\times(0,T')$ with $\Omega'\subset\Omega$ and $T'<T$.
\end{theorem}

\begin{remark}
As a byproduct of Theorem \ref{thm:C1a-gen}, we also get an  regularity result at the initial (time) layer/state for the (standard) $p$-parabolic problems $Hu=0$, by simply applying Theorem \ref{thm:C1a-gen} with $\psi \equiv - \| g \|_{\mathrm{L}^\infty(\partial_p\O)}$, so that  $u > \psi$ in $\O$ and  $Hu = \max\{ Hu  ,\psi- u\} = 0$ in $\O$.
\end{remark}

\subsection{Proofs: arguing from the smallest degenerate scale}\label{sec:HOIG} The main purpose of the function spaces introduced is to be able to keep track, in a quantitative fashion, of the level of degeneracy of $|Du|+s$ in order to state and prove precise regularity estimates in a scale of $C^{1,\alpha}$ spaces. The common theme in the proofs of our main results is that we identify a threshold at which the nondegenerate features of the situation start to dominate, and we then subsequently  divide the proof into two cases, the degenerate and the nondegenerate case; we note that such an idea was exploited in \cite{ALS} for the elliptic $p$-obstacle problem.

To exemplify the philosophy, let us here briefly revisit the interior $C^{1,\alpha}$ estimates for the problem $Hu=0$ with the focus on our intrinsic geometry. Indeed, let   $u$ be a weak solution to $Hu = 0$ in $Q_1$ such that $|u| \leq 1$ in $Q_1$, $u(0,0) = 0$,
and let $\ell(x) \equiv Du(0,0) \cdot x$. Assume that $\mu\equiv|D\ell| + s \leq 1$, and set  $\vp(r) \equiv  \max\{r^\alpha, |D\ell|+s\}$, $\alpha\in (0,\alpha_h)$ where $\alpha_h\in (0,1)$ is the implicit interior regularity threshold just depending on $(n,p,\nu,L)$, and appearing in Theorem \ref{thm:grad -}  below.

First, assume that
$\mu=0$. In this case $\vp(r)=r^\alpha$ and it follows by an argument by contradiction, using blow-up sequences, limiting arguments, and that $\alpha<\alpha_h$,  see Lemma \ref{lem:C1a-deg-int}, that for every $r\in(0,1)$ we have the growth estimate
\begin{equation}\label{eq:C1a-intIntro}
|u | \leq \lambda r^{1+\alpha}\quad\text{in }Q_r(\lambda r^\alpha),
\end{equation}
for some constant $\lambda \equiv \lambda(n,p,\nu,L,\alpha)$, $\lambda > 1$.

Second, assume $\mu>0$. In this case, the underlying philosophy is to rescale the solution such that the rescaled gradient plus the new inhomogeneity constant is of order 1, in order to apply the standard $C^{1,\alpha}$ estimate as stated in Theorem \ref{thm:grad -} below. In this case, we  again first use Lemma \ref{lem:C1a-deg-int} below to conclude, for every $r\in (\varrho,1)$, $\varrho\equiv \mu^{1/\alpha}$, that
\begin{equation}\label{eq:u-l-intIntro}
|u - \ell|  \leq 2\bar\lambda r^{1+\alpha}\quad\text{in }Q_r(\bar\lambda r^\alpha),
\end{equation}
 for some constant  $\bar\lambda\equiv \bar\lambda(n,p,\nu,L,\alpha)$, $\bar\lambda > 1$. Note that if $r\in (\varrho,1)$, then $\vp(r)=r^\alpha$, and $\varrho$ is exactly the threshold below which the nondegenerate features of the situation start to dominate. At the threshold $\varrho$ we have, using \eqref{eq:u-l-intIntro},
 \begin{equation}\label{eq:u-l-intIntro+}
|u - \ell|  \leq c\bar\lambda\mu\varrho\quad\text{in }Q_\varrho(\bar\lambda \mu).
\end{equation}
In particular, introducing the change of variables,
$$(y,\tau)\equiv (x,\mu^{p-2}t),\ \tilde u(y,\tau)\equiv u(x,t)/\mu,$$
we see that $\tilde H\tilde u=0$ in $Q_\varrho(\bar\lambda)$, where $\tilde H$ is a degenerate parabolic operator verifying \eqref{asp} with structure constants $\nu$, $L$ and the inhomogeneity constant $\tilde s \equiv s/\mu$. In particular,  from \eqref{eq:u-l-intIntro+} we deduce that
\begin{equation}\label{eq:u-l-intIntro++}
|\tilde u - \tilde \ell|  \leq c\bar\lambda\varrho\quad\text{in }Q_\varrho(\bar\lambda),
\end{equation}
where now $\tilde \ell(y)\equiv D\tilde u(0,0)\cdot y$ and, importantly,  $\tilde\mu\equiv |D\tilde\ell|+\tilde s=1$. After this rescaling we can apply Theorem \ref{thm:grad -}, see the proof of Lemma \ref{lem:C1a-ndeg-int} for details, and conclude that for every $r\in (0,\rho)$ we have
\begin{equation}\label{eq:u-l-intIntro+++}
|\tilde u - \tilde \ell|  \leq c\bar\lambda r\left(\frac{r}{\varrho} \right)^{\alpha}\quad\text{in }Q_r(\bar\lambda).
\end{equation}
Then, scaling back we are done. Basically all our proofs are structured in this way; based on the threshold $\varrho$  we distinguish between the degenerate case ($|D\ell| + s\approx 0$), and the nondegenerate case ($|D\ell| + s\approx 1$).

\subsection{Organization of the paper} In Section \ref{discuss} we review the literature on obstacle problems for evolutionary $p$-Laplace type operators, we put our results into perspective, and we discuss the sharpness of our results.
Section \ref{prelim} is of preliminary nature and we here introduce the concept of weak solutions, we state the intrinsic Harnack inequality and we provide some technical lemmas. In Section \ref{Grad1} we collect the gradient estimates for the (standard) problem $Hw = \divo  a(Dw) - w_t=0$. In Section \ref{RevisitGrad} we revisit the interior $C^{1,\alpha}$ estimates for the standard problem, discussed in Section \ref{Grad1}, with the focus on intrinsic geometry. In Section \ref{sec:contact}, we establish optimal growth estimates for the obstacle problem at (interior) contact points, and in particular in the intrinsic backward $p$-paraboloids ($P_r^-$) and in the full lower intrinsic cylinders $(Q_r^-)$. In this section we also prove
Theorem \ref{thm:C1a}. In Section \ref{sec:initial} the analysis of Section \ref{sec:contact} is modified to handle the initial layer/state, and Theorem \ref{thm:C1a-ini} and Theorem \ref{thm:C1a-re} are proved. In Section \ref{finalproof} the key estimates established in Sections \ref{RevisitGrad}-\ref{sec:initial} are combined into the proofs of Theorems \ref{thm:C1a-int} and \ref{thm:C1a-gen}.

\section{Review of the literature, our contribution and a discussion of sharpness}\label{discuss}

To put our results into perspective, we note that when $p=2$, then the operator appearing in~\eqref{modello} coincides with the heat operator $\Delta u - u_t$, and hence the obstacle problem in~\eqref{e-obs} becomes an obstacle problem for the heat equation.  In the case of linear uniformly parabolic equations we note that there is an extensive  literature on the existence and regularity of generalized solutions to the obstacle problem in Sobolev spaces, and we refer to  \cite{Friedman75} for details. For obstacle problems involving the heat equation we refer to \cite{CPS}, and for linear parabolic equations with not necessarily constant coefficients we refer to \cite{Bla06b,Bla06a,BDM05,BDM06}.

In contrast, the obstacle problem for the evolutionary $p$-Laplace type operator is considerably less understood in the case
$p\neq 2$. In \cite{KMNobs} the second author, together with Tuomo Kuusi and Giuseppe Mingione, in the case $p>2$ and by showing that solutions have exactly the same degree of regularity as the obstacle,  proved optimal regularity results for obstacle problems involving evolutionary $p$-Laplace type operators with obstacles in intrinsic $C^{0,\alpha}$ spaces, $\alpha\in(0,1]$.  A main ingredient in \cite{KMNobs} was a new intrinsic interpolative geometry allowing for optimal linearization principles via blow-up analysis at contact points. This also opened the way to the proof of a removability theorem for solutions to evolutionary $p$-Laplace type equations. A basic feature of \cite{KMNobs} was that no differentiability in time is assumed on the obstacle; this is in line with the corresponding linear results.

In \cite{ALS} some contributions to the theory of higher order regularity in the obstacle problems for the evolutionary $p$-Laplace operator were established.  More precisely, it is proved that solutions which are Lipschitz in time have $C^{1,1/(p-1)}$ approximation at free boundary points in the spatial variables. In \cite{Sha},   $p$-parabolic free boundary problems were considered,  however the main arguments  in \cite{Sha} for the optimal regularity of solutions are not correct and difficult to repair, see Remark 8 in \cite{ALS}.

As mentioned in the introduction,  this seems to be the first paper where optimal H\"older estimates for the gradient of solutions to $p$-parabolic obstacle problems are established. The key novelty of our work arises in the introduction of the intrinsic geometry that reflects the (non-)degeneracy of the solution at reference points. Incorporating the intrinsic geometry, we prove the optimal growth of the solution from the obstacle, see Theorems \ref{thm:C1a} and \ref{thm:C1a-ini}. Our work also gives new insights to the well-established gradient estimates for the standard problems (i.e., those without obstacles), see Theorem \ref{thm:C1a-int}. As a matter for fact, our H\"older estimates for the gradient at degenerate points hold in much longer time-intervals than the standard ones, at a cost of losing arbitrarily small amount of the H\"older exponent. Our paper can be seen as a novel continuation of the study initiated in \cite{KMNobs} where optimal regularity for obstacle problems was considered in $C^{0,\alpha}$ intrinsic interpolative geometries.

Finally, to establish versions of \cite{KMNobs} and the results of this paper in the singular case, $p<2$, represent interesting projects for future research.

\subsection{Optimal regularity and the intrinsic geometry: a discussion of sharpness}\label{sec:novel} As briefly mentioned, a key feature of obstacle problems involving the $p$-Laplace operator is the improvement of the regularity of the solution at free boundary points. To exemplify,  considering the stationary problem $$\max\{\Delta_p u , \psi - u\} = 0,$$ the solution $u$ is as smooth as the obstacle $\psi$ at free boundary points, i.e., at points on $\partial\{u > \psi\}$, see \cite{Lin, ALS}. In particular, for $C^2$-obstacles $\psi$, solutions $u$ have $C^{1,1}$-expansions at the free boundary points, which of course is optimal. However, in the continuation or non-coincidence set $\{u > \psi\}$, the solutions are (weak) $p$-harmonic functions, hence they are only $C^{1,\alpha}$ regular for some universal $\alpha \in (0,1)$. Thus, the regularity is improved at free boundary points.

The  improvement of regularity at free boundary points in the $p$-obstacle problem is mainly due to the strong minimum principle and the Harnack inequality. The true character of the $p$-Laplace operator only appears when the obstacle is degenerate $(D\psi (x_0)=0)$, since otherwise the linear behavior of the obstacle dictates and makes the equation uniformly elliptic. However, if the obstacle is degenerate, then the solution $u$ is more or less nonnegative, up to an additive constant. Thus, by the strong minimum principle, $u$ cannot deviate too much from above from the obstacle either. Moreover, one can quantify this information sharply by means of the Harnack inequality, see \cite{Lin} (and also \cite{Caff} for the classical obstacle problem for the Laplace operator, $p=2$).

In contrast to the stationary case, it is well-known that for $p>2$, solutions to quasilinear parabolic equations exemplified by the evolutionary $p$-Laplace equation do not satisfy the strong minimum principle, and the Harnack inequality must incorporate  intrinsic geometry. To keep our discussion simple, let us consider the prototype case $$\max\{\Delta_p u - u_t , \psi - u\} = 0.$$ For those who are not familiar with the evolutionary $p$-Laplace equation $\Delta_p u - u_t = 0$, we here present some important nonnegative solutions in the degenerate case $p>2$. For the sake of brevity, we here skip to state the explicit form of the constants appearing in the explicit solutions stated below, instead we refer the interested reader to \cite{DBb}. The following solutions to the degenerate evolutionary $p$-Laplace equation can be found in  \cite{DBb}.
\begin{itemize}
\item The Barenblatt solutions:
\begin{equation}\label{eq:barenblatt}
\Gamma_p(x,t)\equiv \frac{1}{t^{n/\lambda}} \left[ 1 - \gamma_p \left(\frac{|x|}{t^{1/\lambda}}\right)^{\frac{p}{p-1}}\right]_+^{\frac{p-1}{p-2}},\quad t > 0,
\end{equation}
for suitable constants $\gamma_p$ and $\lambda$ depending only on $n$ and $p$.
\item Traveling wave solutions: for every $c>0$, and each $i\in\{1,2,\dots,n\}$,
\begin{equation}\label{eq:travel}
u_c(x,t) =u_c(x_1,...,x_n,t) \equiv K(ct - x_i)_+^{\frac{p-1}{p-2}}
\end{equation}
for some $K \equiv K(c,p)$.
\item Other self-similar solutions: for every $T \in \R$,
\begin{equation}\label{eq:barenblatt2}
u_T (x,t)\equiv c_p \left(\frac{|x|^p}{T-t}\right)^{\frac{1}{p-2}},\quad t < T.
\end{equation}
\end{itemize}

Note that all solutions listed above are nonnegative and do vanish in some part of the interior. This rules out the strong minimum principle from our toolbox. Furthermore, note that nonnegative solutions to  $p$-parabolic problems can be considered as solutions to the corresponding obstacle problems with zero obstacle. This is precisely the point that makes the analysis of the $p$-parabolic obstacle problems so delicate. 

Having these examples in mind, let us explain the optimality of the intrinsic geometry revealed in Theorem \ref{thm:C1a}. Theorem \ref{thm:C1a} states the optimal regularity improvement at free boundary points in the {\it intrinsically scaled backward $p$-paraboloids} (of the form in \eqref{eq:p-para}). The sharpness of the result can be shown with the traveling wave solutions $u_c$ in \eqref{eq:travel} as follows.
\begin{itemize}
\item The origin $(0,0)$ is a free boundary point of $u_c$, and
$$
u_c \leq \lambda r^{1+\frac{1}{p-2}}\quad\text{in } Q_r(\lambda r^{\frac{1}{p-2}}).
$$
This in particular shows that the regularity improvement in the intrinsic cylinders can not go above the order $1+\alpha=1+1/(p-2)$.
\item Nevertheless, taking $\lambda =\lambda(c)$ larger if necessary, a direct computation shows that
$$
u_c \leq \lambda r^{1+\frac{2}{p-2}}\quad\text{in } P_r^-(r^{\frac{2}{p-2}}).
$$
Thus, the regularity is in fact improved in the intrinsic backward $p$-paraboloids to the order $1+\alpha=1+2/(p-2)$. Moreover, such an improvement cannot be extended in the forward time direction.
\end{itemize}

Let us close  with some remarks concerning Theorem \ref{thm:C1a-re}. To distangle the statement in \eqref{eq:C1a-re-asmp}, let $(x_0,t_0)=(0,0)$, $u(0,0)=0$, and assume that  $\bar r=1$, and $\bar\mu=1$,  and let  $\vp(r)=\max\{\omega(r), |D\psi(0,0)| + s \}$. Let $\bar\lambda(r)$ and $\lambda(r)$ be the numerator and denominator, respectively, of the ratio in the left-hand side of \eqref{eq:C1a-re-asmp}, i.e.,
\begin{equation}\label{eq:l-lb}
\begin{aligned}
\bar \lambda(r)  &\equiv \inf\{\lambda > 1: |u|\leq \lambda r\omega(r) \text{ in }Q_r(\omega(r))\},\\
\lambda(r) &\equiv \inf\{\lambda > 1: |u|\leq \lambda r\omega(r) \text{ in }Q_r(\lambda\omega(r))\}.
\end{aligned}
\end{equation}
Note that as $\lambda(r)>1$ we have $Q_r(\lambda(r)\vp(r))\subset Q_r(\vp(r))$. Consequently, we always have $\lambda(r)\leq \bar\lambda(r)$, but in \eqref{eq:C1a-re-asmp} we actually assume that these two quantities are (uniformly) comparable  for all {\it degenerate} scales $r$ ($|D\psi(0,0)| + s\leq\omega(r)$) as measured by $\theta$. Considering the degenerate scales,  and $\omega(r)=r^\alpha$, we have
\begin{align}\label{comp}
 Q_r(\lambda\omega(r))=B_r \times (-\lambda^{2-p}r^{2-\alpha(p-2)}, \lambda^{2-p}r^{2-\alpha(p-2)}),
 \end{align}
 and in this case the assumption in \eqref{eq:C1a-re-asmp} implies that
 \begin{equation}\label{comp+}
 |u| \leq
 \begin{cases} \lambda(r) r\omega(r), & \mbox{in } B_r \times (-(\lambda(r))^{2-p}r^{2-\alpha(p-2)}, (\lambda(r))^{2-p}r^{2-\alpha(p-2)}),\\
 \theta\lambda(r) r\omega(r), & \mbox{in } B_r \times (-r^{2-\alpha(p-2)}, r^{2-\alpha(p-2)}).
 \end{cases}
\end{equation}

As noted in Remark \ref{rem:C1a-re}, \eqref{eq:C1a-re-asmp} is a {\it nontrivial} criterion to extend the optimal improvement of regularity from the intrinsic backward $p$-paraboloids to the entire cylinders. To see this let us  consider the  traveling wave solutions $u_c$ in \eqref{eq:travel} for example. Again let $\lambda(r)$ and $\bar\lambda(r)$ be as in \eqref{eq:l-lb}. Note that both $\lambda(r)$ and $\bar\lambda(r)$ depend on the modulus of continuity $\omega(r)$ for the obstacle $\psi$, and recall that $u_c$ solves the obstacle problem for the  evolutionary $p$-Laplace operator with zero obstacle, as well as that $(x_0,t_0) = (0,0)$ is a degenerate point of $u_c$ on the free boundary, $ct = x_i$. If we choose $\omega(r) \approx r^{2/(p-2)}$, then in the case of $u_c$ we have $\lambda(r) \approx r^{-1/(p-2)}$ and $\bar\lambda (r) \approx r^{-2/(p-2)}$. Hence, $\bar\lambda(r)/\lambda(r) \approx r^{-1/(p-2)}$, so the ratio diverges as $r\to 0^+$, that is, \eqref{eq:C1a-re-asmp} fails to hold. In fact, the ratio $\bar\lambda(r)/\lambda(r)$ diverges whenever we choose $\omega(r) \approx r^\alpha$ with $\alpha > 1/(p-2)$. Nevertheless, the ratio becomes bounded when $\omega(r) \lesssim r^{1/(p-2)}$. Therefore, the criterion in \eqref{eq:C1a-re-asmp} is {\it not} superfluous.



\section{Preliminaries}\label{prelim}


\subsection{Concept of solutions} \label{sec:solution}
If $ U \subset \mathbb R^{n} $ is open and $ 1 \leq q \leq \infty, $ then
by $ W^{1 ,q} ( U ), $ we denote the space of equivalence classes of functions $ f $ with distributional gradient $ Df = ( f_{x_1}, \dots, f_{x_n} ), $ both of which are $ q $-th power integrable on $ U. $ Let
\[ \| f \|_{ W^{1,q} (U)} \equiv  \| f \|_{ \mathrm{L}^q (U)} + \|   | D f |   \|_{ \mathrm{L}^q ( U )}   \]
be the norm in $ W^{1,q} ( U ) $ where $ \| \cdot \|_{\mathrm{L}^q ( U )} $ denotes the usual Lebesgue $q$-norm in $U$.   Given $t_1<t_2$ we denote by $\mathrm{L}^q(t_1,t_2,W^{1,q} ( U ))$ the space of functions such that
for almost every $t$, $t_1\leq t\leq t_2$, the function $x\mapsto u(x,t)$ belongs to
$W^{1,q} ( U )$ and
$$\| u \|_{ \mathrm{L}^q(t_1,t_2,W^{1,q} ( U ))}\equiv \biggl (\iint_{(t_1,t_2)\times U}\bigl (|u(x,t)|^q+|Du(x,t)|^q\bigr )\,\d x\d t\biggr )^{1/q} <\infty.$$ In the following we first describe the concept of weak solutions to
\begin{equation}\label{hom}
Hw = \divo  a(Dw) - w_t=0,
\end{equation}
when the underlying domain considered is not necessarily a cylinder.
\begin{definition}\label{defi1} Let $H$ be as in \eqref{hom} and assume \eqref{asp}. We say that a function $w$
is a weak supersolution (subsolution) to \trif{hom} in an open set $\Xi \Subset \er^{n+1}$ if, whenever
$\Xi'= U \times (t_1,t_2) \Subset \Xi$ with $U\subset\er^n$ and $t_1<t_2$, then
$
w \in \mathrm{L}^p(t_1,t_2;W^{1,p}(U))
$
and
\begin{equation}\label{eq:weak}
\iint_{\Xi'} \left(\langle a(Dw),D\phi\rangle - w \phi_t \right) \,\d x\d t   \geq   (\leq) \ 0,
\end{equation}
for all nonnegative $\phi \in C_0^\infty(\Xi')$. A weak solution is a distributional solution satisfying~\eqref{eq:weak} with equality and without sign restrictions for the test functions. 
\end{definition}

Note that in Definition \ref{defi1} no assumption on the time derivative of $w$ is made. We are now ready to give the definition of solutions to the obstacle problem in \eqref{e-obs}. In the following we assume that the obstacle $\psi$ and boundary value function $g$ are continuous on $\bar \O\subset\mathbb R^{n+1}$ and that $g\geq \psi$ on the parabolic boundary of $\O = \Omega \times (0,T)$.
\begin{definition}\label{defi2} A function
$u$ is a solution to~\eqref{e-obs} if it satisfies the following three properties:
\begin{itemize}
\item[(i)] $u$ is continuous on $\bar \O$, $u \geq \psi$ in $\O$ and $u = g$ on $\partial_{\rm{P}} \O$,
\item[(ii)] $u$ is a weak supersolution in $\O$,
\item[(iii)] $u$ is a weak solution in $\O \cap \{u>\psi\}$.
\end{itemize}
\end{definition}
As for the property (iii), we recall that $u$ is a weak solution in $\O \cap \{u>\psi\}$ means that $u$ is a standard distributional solution in the sense of Definition \ref{defi1} in every space-time cylinder contained in $\O \cap \{u>\psi\}$. We note that a solution to the obstacle problem as in Definition
\ref{defi2} exists by the results in~\cite{KKS}. To be precise, in~\cite{KKS} the boundary values were given by the obstacle itself but it is straightforward to modify the argument in \cite{KKS} to obtain the existence result for general boundary values assuming  $g \geq \psi$ on the parabolic boundary. Moreover, the solution is easily seen to be unique by an ``elliptic" comparison principle for weak solutions, see Lemma \ref{ellittico} below. As outlined in \cite{KMNobs}, there are naturally other ways to obtain existence. An argument arising from potential theory is given in~\cite{LP} and by uniqueness arguments this solution coincides with the solution obtained in~\cite{KKS}. In fact, from~\cite{LP} one finds an argument for an existence result when the obstacle belongs merely to a parabolic Sobolev space. If the obstacle, on the other hand, belongs to parabolic Sobolev space and has time derivative in $\mathrm{L}^2$, then the existence follows from~\cite{AL} and by an approximation argument this approach can be used to obtain the unique solution also in the case when the obstacle is merely continuous; related existence results under regularity assumptions on the obstacle, such as in the existence of $\psi_t$ in suitable Lebesgue spaces, can be found in~\cite{BDM}. Furthermore, an additional approach is given by viscosity solutions in which case the existence is rather easy to obtain. It turns out that a viscosity solution to the obstacle problem is also a so-called $a$-superparabolic function in $\O$, see~\cite{KilLin, KKP, JLM, JJ}, and a continuous weak solution in $\O \cap \{u>\psi\}$. Every bounded superparabolic function is also a weak supersolution by~\cite{KilLin,KKP} and therefore any viscosity solution is a solution in the above sense and hence unique.


\subsection{Technical tools}

Concerning the notion of solution considered above we will several time use the following result (see for example \cite{KilLin} and \cite[Corollary 4.6]{KKP}).
\begin{lemma}[``Elliptic comparison"]\label{ellittico} Let $S \subset \er^{n+1}$ be an open and bounded set and let $T \in\er$. Let $S_T \equiv  S \cap \{t<T\}$. Let $u$ be a weak supersolution in $S_T$ and let $v$ be a weak subsolution in $S_T$. Assume further that $u$ and $v$ are continuous on the closure of $S_T$. If $v \leq u $ on $\partial S_T \setminus \{t=T\}$, then $v \leq u$ in $S_T$.
\end{lemma}

The following Harnack estimate can be retrieved from~\cite{DGV1,DGV2} and~\cite{K}.
\begin{theorem} \label{thm:harnack}
Let $w$ be a nonnegative weak solution to~\eqref{hom} in a space-time cylinder $\O$ in $\mathbb R^{n+1}$. Then there exist positive constants $c$ and $\gamma$, both depending only on $n$, $p$, $\nu$ and $L$, such that whenever
$$
B_{4r}(x_0) \times (t_0- \theta (4r)^p,t_0+ \theta r^p) \subset \O,\quad \text{with}\quad \theta \equiv  \left(\frac{c}{w(x_0,t_0)}\right)^{p-2},
$$
then either
$$
\gamma s r > \min\{1, w(x_0,t_0)\},
$$
or
$$
\gamma^{-1} \sup_{B_r(x_0)} w(\cdot, t_0 - \theta r^p) \leq w(x_0,t_0) \leq \gamma \inf_{B_r(x_0)} w(\cdot,t_0 + \theta r^p).
$$
\end{theorem}

Next we consider linear second order parabolic  equations of the type
\begin{equation}\label{eq:linear eq}Pv\equiv \divo (b(x,t)Dv) -v_t= 0,\end{equation}
where $b$ is a $(n\times n)$-dimensional matrix coefficient satisfying \begin{equation}\label{eq:linear growth} \langle b(x,t)\xi,\xi\rangle \leq \bar L |\xi|^2  , \quad   \langle  b(x,t)\xi,\xi    \rangle \geq \bar \nu |\xi|^2,\end{equation}
whenever $\xi \in \er^n$ and for almost every $(x,t) \in \er^n \times \er$.  Here $0< \bar \nu\leq \bar L$ are fixed parameters.  We will use the following lemma which is an easy consequence of interior $C^{0,\alpha}$ regularity for weak solutions and the strong minimum principle.

\begin{lemma}\label{lemma:Gaussian}Suppose that $v$ is a continuous weak solution to $Pv=0$ in $Q_{2r}^-$ and assume that  $0\leq v \leq 1$ in $Q_{2r}^-$. Then, given $\varepsilon\in (0,1)$,  there exists $\delta \in (0,1)$, depending only on $n$, $\bar \nu$, $\bar L$ and  $\varepsilon$, such that
$$
v(0,0)\leq \delta \quad\Longrightarrow\quad\sup_{Q_r^-} v \leq \varepsilon.
$$
\end{lemma}


\section{Gradient estimates for weak solutions}\label{Grad1}

Here we present the gradient estimates for weak solutions to $Hu = 0$, in a form slightly different compared to the standard formulations, see \cite{DBb}. More specifically, we present them in terms of pointwise approximation by time-independent affine functions. The H\"older exponent $\alpha_h$ appearing below depends on the structure constants, i.e. on $n$, $p$, $\nu$, $L$, and $\alpha_h$ will be fixed throughout the paper.

\begin{theorem}\label{thm:grad -} Let $H$ be as in \eqref{generalH}, with $(a,s)$ as in \eqref{asp}. Let $\ell$ be a time-independent affine function, and let $u$ be a weak solution to $Hu = 0$ in $Q_{4\varrho}(\mu)$, such that
$$
|u - \ell| \leq \lambda\mu \varrho \quad\text{in }Q_{4\varrho}(\lambda\mu)\quad\text{with}\quad  |D\ell| + s\leq \mu,
$$
for some constant $\lambda\geq 1$. Then there exists $c\equiv c(n,p,\nu,L)$, $c \geq 1$, such that
\eqn{grad bound}
$$
|Du| \leq c \lambda \mu\quad\text{in }Q_\varrho(\lambda\mu),
$$
and there exists $\alpha_h\equiv \alpha_h(n,p,\nu,L)$, $\alpha_h\in(0,1)$ such that
\eqn{grad osc}
$$
\osc_{(x_0,t_0) + Q_r(\lambda\mu)}D u \leq c \lambda \mu \left(\frac{r}{\varrho} \right)^{\alpha_h},
$$
whenever $(x_0,t_0) + Q_r(\lambda\mu) \subset Q_\varrho(\lambda\mu)$.
\end{theorem}

To prove Theorem \ref{thm:grad -} we will use the following lemma.

\begin{lemma}[Caccioppoli estimate for differences] \label{lemma:affine caccioppoli}
Let $v$ and $w$ be weak solutions to $Hv=0=Hw$ in $Q_r(\mu)$. Then
\begin{equation}\label{basic energy 1}
\iint_{Q_{r/2}(\mu)} |D(v-w)|^p  \,\d x\d t\leq \frac{c}{r^p}\iint_{Q_r(\mu)} |v-w|^p   \,\d x\d t + c \iint_{Q_r(\mu)}  (\mu + |Dw|+s)^{p}  \,\d x\d t,
\end{equation}
with $c\equiv c(n,p,\nu,L)$.
\end{lemma}
\begin{proof}
 Let $\eta = (v-w)\phi^p$, where $\phi$ is smooth and vanishes outside of $Q_r(\mu)$. Furthermore,  we construct $\phi$ so that $\phi=1$ on $Q_{r/2}(\mu)$,  $|D\phi| \leq 4/r$ and $|\phi_t| \leq 8 \mu^{p-2}/r^2$ on $Q_r(\mu)$, and $0\leq \phi \leq 1$ in $Q_r(\mu)$.  As $Hv=0=Hw$ weakly in $Q_r(\mu)$  we have, formally,
 \begin{align*}
0  =& \iint_{Q_r(\mu)} \left( \langle a(Dv)-a(Dw) , D \eta \rangle -(v-w) \partial_t \eta  \right)   \,\d x\d t=  I_1 + I_2 - I_3,
\end{align*}
where
\begin{align*}
I_1\equiv &\iint_{Q_r(\mu)} \langle a(Dv)-a(Dw) , D(v-w) \rangle \phi^p   \,\d x\d t, \\
I_2\equiv &p \iint_{Q_r(\mu)} \langle a(Dv)-a(Dw) , D\phi\rangle (v-w) \phi^{p-1}   \,\d x\d t, \\
I_3\equiv &\frac12 \iint_{Q_r(\mu)} (v-w)^2 \partial_t \phi^p   \,\d x\d t,
\end{align*}
after an integation by parts in time. The integration by parts argument can be made rigorous using Steklov averages but we here omit the routine details. Using  the assumptions in \eqref{asp} on $a(\cdot)$,
 $$
\begin{aligned}
I_1 & \geq  \frac1c  \iint_{Q_r(\mu)}|D(v-w)|^p  \phi^p   \,\d x\d t.
\end{aligned}
$$
Using further structure of $a(\cdot)$,
$$
\begin{aligned}
|I_2| \leq & c \iint_{Q_r(\mu)}\left(|Dv| + |Dw| + s \right)^{p-2} |D(v-w)| \phi^{p-1} |v-w| |D\phi|   \,\d x\d t\\
 \leq&
 c \iint_{Q_r(\mu)} (|Dw|+s)^{p-2}  |D(v-w)| \phi^{p-1} |v-w| |D\phi|  \,\d x\d t \\
 &+c \iint_{Q_r(\mu)} |D(v-w)|^{p-1} \phi^{p-1} |v-w| |D\phi|  \,\d x\d t.
\end{aligned}
$$
Using Young's inequality repeatedly,
\begin{align*}
&(|Dw|+s)^{p-2}  |D(v-w)| \phi^{p-1} |v-w| |D\phi|\\
&\leq (|Dw|+s)^{p-2}  |D(v-w)| \phi |v-w| |D\phi|\\
&\leq c(p)\bigl ((|Dw|+s)^{p}+ |D(v-w)|^{p/2} \phi^{p/2} |v-w|^{p/2} |D\phi|^{p/2}\bigr )\\
&\leq c(p)\bigl ((|Dw|+s)^{p}+ \varepsilon^p|D(v-w)|^{p} \phi^{p} + \varepsilon^{-p}|v-w|^{p} |D\phi|^{p}\bigr ),
\end{align*}
where $\varepsilon\in (0,1)$ is a degree of freedom. Combining this with one more application of  Young's inequality, in the second term in the estimate of $|I_2|$ above, we deduce
\begin{align*}
|I_2|
 \leq &
c\varepsilon^p \iint_{Q_r(\mu)} |D(v-w)|^p  \phi^p \,\d x\d t+ c\tilde c(\varepsilon)  \iint_{Q_r(\mu)}(|Dw|+s)^p   \,\d x\d t\\
& +c r^{-p} \iint_{Q_r(\mu)}|v-w|^p    \,\d x\d t,
\end{align*}
with $c\equiv c(n,p,\nu,L)$. Again, using Young's inequality,
 $$
\begin{aligned}
|I_3| & \leq \frac12 \iint_{Q_r(\mu)} (v-w)^2 |\partial_t \phi^p|   \,\d x\d t \leq  c \mu^{p-2}  r^{-2}  \iint_{Q_r(\mu)}(v-w)^2   \,\d x\d t
  \\ \nonumber
 & \leq  cr^{-p}  \iint_{Q_r(\mu)} |v-w|^p   \,\d x\d t  + c \mu^p |Q_r(\mu)|,
\end{aligned}
$$
with $c\equiv c(p)$.  Collecting estimates leads to the statement of the lemma.
\end{proof}

We are now ready to prove Theorem \ref{thm:grad -}.

\begin{proof}[Proof of Theorem \ref{thm:grad -}]
Let in the following $c$ be a positive constant, depending only on $n$, $p$, $\nu$, and $L$, which is allowed change upon each occurrence. Applying Lemma~\ref{lemma:affine caccioppoli} to the pair $u$ and $\ell$ in $Q_{4\varrho}(\lambda\mu)$, recalling that $|D\ell| + s\leq\mu$, $\lambda\geq 1$, we obtain
$$
\begin{aligned}
\bariint_{Q_{2\varrho}(\lambda\mu)} |Du - D\ell |^p   \,\d x\d t & \leq  \frac{c}{\varrho^p} \bariint_{Q_{4\varrho}(\lambda\mu)} |u-\ell|^p   \,\d x\d t + c(\lambda\mu)^p  \leq  c(\lambda\mu)^p.
\end{aligned}
$$
Now, using the Minkowski and H\"older inequalities we deduce
$$
\bariint_{ Q_{2\varrho}(\lambda\mu) } |Du |^{p-1}   \,\d x\d t  \leq c(\lambda\mu)^{p-1} + |D\ell|^{p-1} \leq c (\lambda\mu)^{p-1}.
$$
Employing the local boundedness for the gradient, e.g., Theorem 5.1 in Ch. VIII in \cite{DBb}, we obtain
$$
|Du| \leq c \lambda\mu\quad\text{in }Q_{3\varrho/2}(\lambda\mu),
$$
which proves~\eqref{grad bound}. We can then use the well known H\"older estimate for the gradient, e.g., Theorem 1.1 in Ch. IX in \cite{DBb}, to complete the proof.
\end{proof}



\section{Interior gradient estimates revisited}\label{RevisitGrad}

In this section we revisit the interior $C^{1,\alpha}$ estimates for the (standard) problem $Hu=0$ with the focus on intrinsic geometry. In the following we let  $\alpha_h$ denote the H\"older exponent for the gradient, appearing in Theorem \ref{thm:grad -}, and we let $\alpha$ be any exponent satisfying
\begin{equation}\label{eq:alpha}
0< \alpha < \alpha_h.
\end{equation}
The purpose of this section is to prove the following proposition.

\begin{proposition}\label{prop:C1a-int} Let $H$ be as in \eqref{generalH}, with $(a,s)$ as in \eqref{asp}, and let $u$ be a weak solution to $Hu = 0$ in $Q_1$ such that $|u| \leq 1$ in $Q_1$, $u(0,0) = 0$, and $|Du(0,0)| + s \leq 1$. Let $\alpha \in (0,\alpha_h)$ be given, set $\vp(r) \equiv  \max\{r^\alpha, |Du(0,0)|+s\}$, and let $\ell \equiv Du(0,0) \cdot x$. Then for every $r\in(0,1)$,
\begin{equation}\label{eq:C1a-int}
|u - \ell | \leq \lambda r^{1+\alpha}\quad\text{in }Q_r(\lambda\vp(r)),
\end{equation}
for some constant $\lambda \equiv \lambda(n,p,\nu,L,\alpha)$, $\lambda > 1$.
\end{proposition}


\subsection{The degenerate case}\label{sec:deg-int} We first treat the degenerate case, $|Du(0,0)| + s\approx 0$.

\begin{lemma}\label{lem:C1a-deg-int}
Under the assumptions as in the statement of Proposition \ref{prop:C1a-int}, there exists $\lambda\equiv \lambda(n,p,\nu,L,\alpha)$, $\lambda>1$, such that if $|Du(0,0)| + s < r^\alpha < 1$, then
\begin{equation}\label{eq:C1a-deg-int}
| u | \leq \lambda r^{1+\alpha}\quad\text{in }Q_r(\lambda r^\alpha).
\end{equation}
\end{lemma}

\begin{proof} We are going to prove that  if $H$, $s$, $\alpha$, and $u$  are as in the setting of the lemma, then there exist  $\bar\lambda \gg 1$ and $\bar r$, $0<\bar r\ll 1$, both determined {\it a priori} by $n$, $p$, $\nu$, $L$ and $\alpha$ only, such that if for some $\lambda > \bar\lambda$, and for some  $r<\bar r$, $|Du(0,0)| + s < r^\alpha < 1$,  we have that
\begin{equation}\label{eq:C1a-deg-int-asmp}
|u| < \lambda \varrho^{1+\alpha}\quad\text{in }Q_\varrho(\lambda\varrho^\alpha),
\end{equation}
for all $\varrho \in (r,1)$, then the strict inequality in \eqref{eq:C1a-deg-int-asmp} continues to hold at $\varrho = r$. Once such constants are found, the proof can be easily completed by our initial assumption that $|u|\leq 1$ in $Q_1$, and  a standard continuity argument, see Remark \ref{MethodCont} below.

With this goal we instead assume, by way of contradiction, that there exist, for each $j=1,2,\dots$, $H_j$, $s_j$, $u_j$ as in the setting of the lemma,  some $r_j\to 0$ and $\lambda_j\to\infty$, such that $|Du_j(0,0)| + s_j < r_j^\alpha$, and such that
\begin{equation}\label{eq:C1a-deg-int-asmp-re}
|u_j| < \lambda_j \varrho^{1+\alpha}\quad\text{in }Q_\varrho(\lambda_j\varrho^\alpha),
\end{equation}
for all $\varrho\in(r_j,1)$, but
\begin{equation}\label{eq:C1a-deg-int-false}
|u_j(x_j,t_j)| \geq \lambda_j r_j^{1+\alpha},
\end{equation}
for some $(x_j,t_j) \in Q_{r_j}(\lambda_jr_j^\alpha)$. Using Theorem \ref{thm:grad -}, with $\ell = 0$, $\lambda = \lambda_j$, and
$\mu=\varrho^\alpha$, we see that \eqref{eq:C1a-deg-int-asmp-re} implies that
\begin{equation}\label{eq:Duj}
|Du_j(\cdot,t) - Du_j(0,t)| \leq c\lambda_j\varrho^\alpha \left(\frac{r_j}{\varrho}\right)^{\alpha_h}\quad\text{in }B_{2r_j},
\end{equation}
whenever $\varrho \in (8r_j,1)$ and $|t|\leq (2r_j)^2(\lambda_j\varrho^\alpha)^{2-p}$.

Define
\begin{equation}\label{eq:tuj}
\tilde u_j (y,\tau) \equiv \frac{u_j(r_j y,r_j^2\theta_j^{2-p}\tau)}{r_j\theta_j}\quad\text{with}\quad \theta_j \equiv \lambda_j r_j^\alpha.
\end{equation}
Then $\tilde u_j$ is a weak solution to
$$
\tilde H_j \tilde u_j = 0 \quad \text{in }Q_2(2^\alpha),
$$
where $\tilde H_j$ verifies \eqref{asp} with $\nu$, $L$ and $s \equiv s_j/\theta_j \to 0$. By \eqref{eq:C1a-deg-int-asmp-re} (with $\varrho = 2r_j$), we have
$$
|\tilde u_j| \leq 2^{1+\alpha}\quad\text{in }Q_2(2^\alpha).
$$
Hence, by the interior regularity theory, see \cite{DBb}, both $\{\tilde u_j\}$ and $\{D \tilde u_j\}$ are uniformly H\"older continuous in $\overline{Q_1}$. Then by the Arzela-Ascoli theorem, $\tilde u_j\to \tilde u$ and $D\tilde u_j\to D \tilde u$ uniformly in $\overline{Q_1}$ along a subsequence, for some continuous function $\tilde u : \overline{Q_1} \to \R$ having continuous spatial derivatives.

Writing $\tilde H_jw\equiv \ddiv \tilde a_j(Dw) - w_\tau$, it follows from \eqref{asp} that $\{\tilde a_j\}$ is also a locally Lipschitz continuous map in $\R^n$. Hence, extracting a further subsequence along which $\tilde u_j\to \tilde u$ and $D\tilde u_j\to D \tilde u$, we have $\tilde a_j \to \tilde a$ locally uniformly in $\R^n$, for some locally Lipschitz vector-field $\tilde a:\R^n\to\R^n$. In particular,
$$
\tilde a_j(D\tilde u_j) \to \tilde a (D\tilde u) \quad\text{uniformly in }\overline{Q_1}.
$$
Hence, we must have
$$
\tilde H\tilde u = 0\quad\text{in }Q_1,
$$
in the weak sense, where $\tilde Hw\equiv \ddiv \tilde a(Dw) - w_\tau$.

Rewriting \eqref{eq:Duj} in terms of $\tilde u_j$,
\begin{equation}\label{eq:u-C1a-deg-int-re}
|D\tilde u_j (\cdot,\tau) - D\tilde u_j(0,\tau)| \leq c r_j^{\alpha_h-\alpha}\quad\text{in }\overline{B_1},
\end{equation}
whenever $|\tau| \leq 1$. Since $\alpha < \alpha_h$ and $r_j\to 0$, the uniform convergence of $D\tilde u_j\to D\tilde u$ (along some subsequence) implies that
$$
D\tilde u (\cdot, \tau) \equiv D\tilde u(0,\tau)\quad\text{in }\overline{B_1},
$$
whenever $|\tau| \leq 1$. Thus, for every test function $\vp\in C_0^\infty(Q_1)$,
$$
\iint_{Q_1} \tilde u \vp_\tau \,\d y \d\tau = \iint_{Q_1} \tilde a(D\tilde u)\cdot D\vp\,\d y \d\tau = -\iint_{Q_1}\ddiv \tilde a(D\tilde u) \vp \,\d y \d\tau = 0,
$$
which implies that $\tilde u(y,\tau) \equiv \tilde u(y,0)$ for every $(y,\tau)\in Q_1$. However, by \eqref{eq:u-C1a-deg-int-re}, $\tilde u_j(0,0) = 0$ (by assumption) and $|D\tilde u_j(0,0)| \leq r_j^\alpha \to 0$, we have that $\tilde u(\cdot,0) \equiv 0$ on $\overline{B_1}$,  and hence
$$
\tilde u \equiv 0\quad\text{in }\overline{Q_1}.
$$
Nevertheless, rephrasing \eqref{eq:C1a-deg-int-false} in terms of $\tilde u_j$ and passing to the limit, there must exist a point $(y_0,\tau_0)\in \overline{Q_1}$ such that  $|\tilde u (y_0,\tau_0)| \geq 1$, and this gives a contradiction.\end{proof}

\begin{remark} \label{MethodCont} In the proof of Lemma \ref{lem:C1a-deg-int} we refer to a (standard) continuity argument. To outline the argument, assume the hypotheses of Lemma \ref{lem:C1a-deg-int} and that it is true that if
\begin{equation}\label{keyasump}
|u| < \lambda \varrho^{1+\alpha}\quad\text{ in }Q_\varrho(\lambda\varrho^\alpha),
\end{equation}
for all $\varrho \in (r,1)$,  then also
\begin{equation}\label{keyasump+} |u| < \lambda r^{1+\alpha}\quad\text {in }Q_r(\lambda r^\alpha).
\end{equation}
Let $\delta\in (0,1)$ be a given number, and let
\begin{equation}\label{keyasump++}
I \equiv \{ r \in [\delta,1]:\ |u| < \lambda \varrho^{1+\alpha}\quad\text{ in }Q_\varrho(\lambda \varrho^\alpha),\quad \text{ for all } \varrho \in [r,1]\}.
\end{equation}
Observe that $1 \in I$, because we assume $|u| \leq 1$ in $Q_1$, hence $I$ is nonempty. If $r \in I$, then by the continuity of $u$, $(r-\e, r+\e) \subset I$, for some $\e > 0$ (which can be very small, and dependent on $u$ and $r$), and hence $I$ is open. We also claim that $I$ is a closed. Indeed, suppose $\{r_k\} \subset I$ is such that $r_k \to r_0$. If $r_0 > \inf_k r_k$, then by definition, $r_0 \in I$. Otherwise, if $r_0 = \inf_k r_k$, then by continuity of $|u|$, $|u| < \lambda \varrho^{1+\alpha}$ in $Q_\varrho$, $\forall \varrho \in (r_0, 1]$. This is true, because for any $\hat r_0 > r_0$, we can find some $r_k < \hat r_0$, so $|u| < \lambda \varrho^{1+\alpha}$ in $Q_\varrho$,  $\forall \varrho \in [r_k,1] \supset [\hat r_0,1]$. Now the implication in \eqref{keyasump}-\eqref{keyasump+} applies and $r_0 \in I$. Thus, $I$ is a closed set. Having proved that $I$ is a nonempty, open and closed set, we can conclude that $I \equiv [\delta, 1]$.
\end{remark}

\subsection{The nondegenerate case}\label{sec:ndeg-int} We here consider the nondegenerate case, $|Du(0,0)| + s \approx 1$.

\begin{lemma}\label{lem:C1a-ndeg-int} Assume, in addition to the assumptions in Proposition \ref{prop:C1a-int}, that $\mu \equiv  |Du(0,0)| + s > 0$, and that there exists
$\varrho \in (0,1)$, $\varrho^\alpha\leq \mu$, such that
$$
| u - \ell | \leq \lambda \mu \varrho\quad\text{in }Q_\varrho(\lambda\mu),
$$
for some $\lambda > 1$. Then there exists a constant $c\equiv c(n,p,\nu,L)$, $c> 1$, such that if $r\in(0,\varrho)$, then
$$
|u - \ell | \leq  c\lambda \mu r\left(\frac{r}{\varrho}\right)^\alpha \quad\text{in }Q_r(\lambda \mu).
$$
\end{lemma}

\begin{proof}
By Theorem \ref{thm:grad -} and the fact that $D\ell = Du(0,0)$, we obtain
\begin{equation}\label{eq:Du-Dl}
|D u - D\ell| \leq c_0 \lambda \mu \left(\frac{r}{\varrho}\right)^\alpha,\quad\forall r\in (0,\varrho/4).
\end{equation}
Also recalling that $u(0,0) = \ell(0) = 0$, that $\ell$ is time-independent, we deduce
$$
\sup_{Q_r(\lambda\mu)} |u - \ell | \leq \sup_{-r^2(\lambda\mu)^{2-p} < t < r^2(\lambda\mu)^{2-p} }\left| \means{B_r} (u (x,t) - u (x,0))\,\d x \right| + 2c_0\lambda  \mu r \left(\frac{r}{\varrho}\right)^\alpha,
$$
whenever $4r < \varrho$. To estimate the first term on the right-hand side in the last display, we use the weak formulation of $Hu = 0$ in $Q_\varrho(\lambda\mu)$ and the divergence theorem. That is,
$$
\begin{aligned}
\means{B_r} (u (x,t) - u (x,0))\,\d x &= \means{B_r}\,\d x \int_0^t\frac{\partial u}{\partial t} (x,\tau) \,\d\tau \\
& = \int_0^t \d\tau \means{B_r} \ddiv a(Du)\,\d x = \frac{n}{r} \int_0^t\,\d\tau \means{\partial B_r}   a(Du)\cdot \vec\nu\,\d\sigma_x,
\end{aligned}
$$
where $\vec\nu$ denotes the outward unit normal to $\partial B_r$. Now we can make use of the structure condition in \eqref{asp}, along with \eqref{eq:Du-Dl} and that $\mu \equiv |D\ell| + s \leq 1$, to deduce that
$$
\begin{aligned}
\left| \means{\partial B_r} a(Du)\cdot \vec\nu \, \d\sigma_x\right| &= \left| \means{\partial B_r} (a(Du) - a(D\ell))\cdot\vec\nu \, \d\sigma_x \right| \\
& \leq \means{\partial B_r} |Du - D\ell| \, \d\sigma_x \int_0^1 |\partial_z a (\eta Du + (1-\eta) D\ell)|\, \d\eta  \leq c_1c_0 L (\lambda \mu)^{p-1} \left(\frac{r}{\varrho}\right)^\alpha,
\end{aligned}
$$
for some $c_1\equiv c_1(p)$. Thus, combining the last three displays, we observe that
$$
\sup_{Q_r(\lambda\mu)} |u - \ell| \leq c_2 c_0 \lambda \mu r \left(\frac{r}{\varrho} \right)^\alpha,
$$
with $c_2 \equiv \max\{ nc_1, 2\}$. Recall that $r$ was an arbitrary number in $(0,\varrho/4)$. Replacing $c_2$ with $4^{1+\alpha} c_2$, we can extend the lasy estimate to all $r\in(0,\varrho)$, as desired.
\end{proof}


\subsection{Proof of Proposition \ref{prop:C1a-int}}\label{sec:pf-C1a-int} The idea underlying the proof is that once we rescale our solution such that the rescaled gradient is of order 1, then we can apply the standard $C^{1,\alpha}$ estimate as formulated in Theorem \ref{thm:grad -}. The only case in which such a rescaling is impossible is when the gradient vanishes, and in that case we have, by Lemma \ref{lem:C1a-deg-int}, the intrinsically scaled $C^{1,\alpha}$ estimate all the way to the origin. A similar idea also appears in the setting of elliptic problems, see \cite{ALS}.

To start the proof, we first use Lemma \ref{lem:C1a-deg-int}. Indeed, if $\mu\equiv |Du(0,0)| + s = 0$, then \eqref{eq:C1a-deg-int} gives \eqref{eq:C1a-int} directly. Thus, it suffices to consider the case $\mu > 0$. We choose $\ell \equiv Du(0,0)\cdot x$. Then by \eqref{eq:C1a-deg-int},
\begin{equation}\label{eq:u-l-int}
|u - \ell| \leq \osc_{Q_r^-(\bar\lambda r^\alpha)} u + \osc_{B_r} \ell \leq (\bar\lambda+1) r^{1+\alpha}\quad\text{in }Q_r(\bar\lambda r^\alpha),
\end{equation}
for every $r\in (\mu^{1/\alpha},1)$, where $\bar\lambda\equiv \bar\lambda(n,p,\nu,L,\alpha)$. Next, we invoke Lemma \ref{lem:C1a-ndeg-int} (with $\lambda = \bar\lambda r^\alpha$ and $\varrho = \mu^{1/\alpha}$ there), which implies that
\begin{equation}\label{eq:u-l-int-re}
|u - \ell| \leq c_1c_0 \bar\lambda r^{1+\alpha}\quad\text{in }Q_r^-(\bar\lambda\mu),
\end{equation}
for every $r\in (0,\mu^{1/\alpha})$, with $c_1 \equiv c_1(n,p,\nu,L)$. By combining these two estimates we get the desired conclusion.


\section{Analysis at contact points}\label{sec:contact}

This section is devoted to the study of optimal growth estimates at contact points. First, we establish the optimal estimate, up to the critical exponent $2/(p-2)$, in the intrinsic backward $p$-paraboloids ($P_r^-$).

\begin{proposition}\label{prop:C1a-para}
Let $H$ be as in \eqref{generalH}, with $(a,s)$ as in \eqref{asp}, and let $\omega$ verify \eqref{eq:omega prop} with $\alpha\in (0,{2}/{(p-2)}]$. Let $\psi$ be a continuous obstacle, let $\ell$ be a time-independent affine function, and let $u$ be a solution to $\max\{ Hu, \psi- u\} = 0$ in $Q_1$ such that $u(0,0) = \psi(0,0) = 0$, and $|D\ell| + s \leq 1$. Let $\vp(r) \equiv \max\{ \omega(r), |D\ell| + s \}$, and assume that
\begin{equation}\label{eq:p-C1a}
|\psi - \ell| \leq r\omega(r) \quad\text{in }Q_r(\vp(r)),
\end{equation}
for every $r\in(0,1)$. Then there exists $\lambda \equiv \lambda(n,p,\nu,L,A,\alpha)$, $\lambda > 1$, such that
\begin{equation}\label{eq:u-C1a-para}
|u | \leq \lambda r\vp(r)\quad\text{in }P_r^-(\vp(r)),
\end{equation}
for every $r\in(0,1)$.
\end{proposition}

Second, we establish the optimal estimate, up to the critical exponent $2/(p-2)$,  in the full lower intrinsic cylinders $(Q_r^-)$ but under an additional intrinsic assumption.

\begin{proposition}\label{prop:C1a-relower} Assume, in addition to the assumptions in Proposition \ref{prop:C1a-para}, that $|u|\leq 1$ in $Q_1$, and that there exists a  constant $\theta\geq 1$ such that if $|D\ell| + s < \omega(r)< 1$, then
\begin{equation}\label{eq:criterion}
\frac{ \inf\{ \lambda > 1:  |u | \leq \lambda r\omega(r)\text{ in }Q_r(\lambda\omega(r))\} }{\inf\{ \lambda > 1: |u| \leq \lambda r\omega(r) \text{ in }Q_r(\omega(r))\} } \leq \theta.
\end{equation}
Then there exists $\bar\lambda\equiv\bar\lambda(n,p,\nu,L,A,\alpha,\theta)$ such that
\begin{equation}\label{eq:u-C1a-re}
|u - \ell | \leq \bar\lambda r\omega(r)\quad\text{in }Q_r^-(\bar\lambda\vp(r)),
\end{equation}
for every $r\in(0,1)$.
\end{proposition}

Third, we establish the optimal estimate, this time up to the  universal exponent $\alpha_h$, in the full lower intrinsic cylinders $(Q_r^-)$ but (essentially) without additional assumptions.

\begin{proposition}\label{prop:C1alower} Assume, in addition to the assumptions in  Proposition \ref{prop:C1a-para}, that $\alpha\in(0,\alpha_h)$, where $\alpha_h$ is as in Theorem \ref{thm:grad -}, and that $|u| \leq 1$ in $Q_1$. Then there exists $\bar\lambda\equiv\bar\lambda(n,p,\nu,L,A,\alpha)$ such that
\begin{equation}\label{eq:u-C1alower}
|u - \ell | \leq \bar\lambda r\omega(r) \quad\text{in }Q_r^-(\bar\lambda\vp(r)),
\end{equation}
for every $r\in(0,1)$.
\end{proposition}

We emphasize that the conclusions of Propositions \ref{prop:C1a-para}-\ref{prop:C1alower} are stated in the sets $$\mbox{$P_r^-(\vp(r))$, $Q_r^-(\bar\lambda\vp(r))$, and  $Q_r^-(\bar\lambda\vp(r))$,}$$ respectively. Furthermore, Proposition \ref{prop:C1a-para} and Proposition \ref{prop:C1a-relower}
are valid for $\alpha\in (0,{2}/{(p-2)}]$, while  Proposition \ref{prop:C1alower} is valid for $\alpha\in (0,\alpha_h)$. The proofs of Propositions \ref{prop:C1a-para}-\ref{prop:C1alower} will be split into several steps which we below present in subsections and we again divide the argument into the degenerate case ($|D\ell| + s \approx 0$), and the nondegenerate case ($|D\ell| + s \approx 1$).

For the degenerate scales $|D\ell| + s < \omega(r) < 1$, we actually prove more in the sense that the conclusions in Proposition \ref{prop:C1a-relower} and Proposition \ref{prop:C1alower} remain true
with $Q_r^-(\bar\lambda\vp(r))$ replaced by $Q_r(\bar\lambda\vp(r))$.  In particular, in this sense, for degenerate scales stronger versions of Proposition \ref{prop:C1a-relower} and Proposition \ref{prop:C1alower} are valid. The restriction to the lower intrinsic cylinders $Q_r^-(\bar\lambda\vp(r))$ in the statement of these propositions appears in our proofs for nondegenerate scales
$\omega(r)\leq |D\ell| + s$, see Lemma \ref{lem:C1a-ndeg}. In fact, in the context of Proposition \ref{prop:C1alower} we in Lemma \ref{lem:C1a-ndeg} prove,  by a linearization of the problem, and by making use of the strong minimum principle for linear second order parabolic equations, that
Proposition \ref{prop:C1alower}  holds for all $\alpha\in (0,{2}/{(p-2)}]$ for nondegenerate scales. Now the restriction to $Q_r^-(\bar\lambda\vp(r))$ arises in the use of the (strong) minimum principle. The proof of Proposition \ref{prop:C1a-relower} for nondegenerate scales is identical to
the one for Proposition \ref{prop:C1alower}, and this explains our restriction to $Q_r^-(\bar\lambda\vp(r))$ in the statements of these propositions.

As we will see, once we have proved Proposition \ref{prop:C1a-ini} below, which concerns the regularity for our problem at the initial (time) layer/state, then by
Proposition \ref{prop:C1a-ini}, Proposition \ref{prop:C1a-relower} and Proposition \ref{prop:C1alower}, we can conclude that the latter two propositions also hold in
the full intrinsic cylinders $Q_r(\bar\lambda\vp(r))$, we refer to Subsection \ref{extension} for details.


\subsection{The degenerate case} We here study the degenerate case, $|D\ell| + s \approx 0$. We first establish the optimal estimate in the intrinsic backward $p$-paraboloids, for all exponents $\alpha \leq {2}/{(p-2)}$.

\begin{lemma}\label{lem:C1a-deg-para}
Under the assumptions stated in Proposition \ref{prop:C1a-para}, for every $\alpha \in (0,{2}/{(p-2)}]$ there exists $\lambda\equiv \lambda(n,p,\nu,L,A,\alpha)$, $\lambda\geq 1$, such that if $|D\ell| + s < \omega(r) < 1$, then
$$
| u | \leq \lambda r\omega(r)\quad\text{in }P_r^-(\omega(r)).
$$
\end{lemma}

\begin{proof}
It is enough to handle the case $\omega^{-1}(|D\ell| + s) \leq \frac{1}{8}$, and hence we fix $r\leq \frac{1}{8}$ such that $|D\ell| + s \leq \omega(r)$. It then follows from \eqref{eq:p-C1a}, and that $\psi(0,0) = 0$,  that
\begin{equation}\label{eq:p-C1a-re}
|\psi| \leq 16r\omega(8r) \quad\text{in }Q_{8r}(\omega(8r)).
\end{equation}
As in \cite{KMNobs} we now construct a comparison map. Consider the following auxiliary problem,
$$
\begin{dcases}
Hv = 0 &\text{in }Q_{8r}(\omega(8r)),\\
v = u + 16 r \omega(8r) & \text{on }\partial_{\rm P} Q_{8r}(\omega(8r)).
\end{dcases}
$$
Since $u\geq \psi$ in $Q_1$ and $Q_{8r}(\omega(8r)) \subset Q_1$ for any $r \leq 1$, it follows from \eqref{eq:p-C1a-re} that
$$
u \geq \psi \geq - 16r\omega(8r)\quad\text{in }Q_{8r}(c\bar\lambda\mu).
$$
Since $Hu\leq 0$ in $Q_1$ in the weak sense, we deduce from the comparison principle (Lemma \ref{ellittico}) that
\begin{equation}\label{eq:v-u-deg}
0 \leq v \leq u + 16r\omega(8r)\quad\text{in }Q_{8r}(\omega(8r)).
\end{equation}

On the other hand, by \eqref{eq:p-C1a-re} and \eqref{eq:v-u-deg},
$$
v \geq 0  \geq u - 16 r \omega(8r)\quad\text{on }\{ u = \psi\} \cap Q_{8r}(\omega(8r)),
$$
and hence, in particular,
$$
v + 16r\omega(8r) \geq u \quad\text{on }(\partial\{u > \psi\} \cap Q_{8r}(\omega(8r)))\setminus (B_{8r}\times \{  (8r)^2\omega(8r)^{2-p}\}).
$$
Since $H v = 0 = H u$ in $\{ u >\psi\}\cap Q_{8r}(\omega(8r))$, $v > u$ on $\partial_{\rm P} Q_{8r}(\omega(8r))$ and $\{ u >\psi\}$ is an open set, we can invoke the comparison principle (Lemma \ref{ellittico}) again, to deduce that
$$
v + 16r\omega(8r)  \geq u \quad\text{in }\{u > \psi\}\cap Q_{8r}(\omega(8r)).
$$
Collecting the last three displays, we arrive at
\begin{equation}\label{eq:v-u-deg-re}
v \geq u - 16r\omega(8r)\quad\text{in }Q_{8r}(\omega(8r)).
\end{equation}

Now since $u(0,0) = \psi(0,0) = 0$, \eqref{eq:v-u-deg} implies
\begin{equation}\label{eq:v-0-deg}
v(0,0) \leq 16r\omega(8r) \leq 2^{7+3\alpha} A r\omega(r),
\end{equation}
where the last inequality is due to \eqref{eq:omega prop}. We claim that
\begin{equation}\label{eq:v-claim}
v(0,t) \leq 2^{7+3\alpha} A\gamma r\omega(r),\quad\forall t\in [-r^2\omega(r)^{2-p},0],
\end{equation}
where $\gamma \equiv \gamma(n,p,\nu,L)$, $\gamma > 1$, comes from the intrinsic backward Harnack inequality, see Theorem \ref{thm:harnack}.

To prove \eqref{eq:v-claim}, it is sufficient to consider the case $v(0,0) > 0$, since otherwise one can first repeat the following argument with $v_\e \equiv v + \e$, and then send $\e\to 0^+$. Thus, let us assume that $v(0,0) > 0$. Then by the backward intrinsic Harnack inequality, see Theorem \ref{thm:harnack},
\begin{equation}\label{eq:forward-harnack}
\sup_{B_\varrho} v (\cdot, - \theta \varrho^p) \leq \gamma v(0,0) \leq 16 \gamma r\omega(8r),
\end{equation}
provided that
\begin{equation}\label{eq:forward-harnack-re}
B_{4\varrho} \times (-\theta(4\varrho)^p, \theta \varrho^p)\quad\text{with}\quad \theta \equiv \left(\frac{c}{v(0,0)}\right)^{p-2},
\end{equation}
is contained in $Q_{8r}(\omega(8r))$; here $c \equiv c(n,p,\nu,L)$, $c>1$, is  the waiting time constant for the intrinsic Harnack inequality. We remark that we may increase the waiting time constant $c$ such that $c > 2^{7+3\alpha} A$, by increasing the Harnack constant $\gamma$, see Ch. 5, paragraph 2.4 in \cite{DGV3}; clearly both $\gamma$ and $c$ will now depend further on $A$ and $\alpha$. Then by \eqref{eq:v-0-deg}, we have
$$
\theta r^p \geq \left( \frac{c}{2^{7+3\alpha} A} \right)^{p-2} r^2 (\omega(r))^{2-p} > r^2 (\omega(r))^{2-p}.
$$
Thus, given any $t \in [-r^2\omega(r)^{2-p},0)$, we can choose $\varrho \equiv\varrho(t) \in (0,r)$ such that $-t = \theta \varrho^p$. With such a choice of $\varrho = \varrho(t)$, the set inclusion in \eqref{eq:forward-harnack-re} is satisfied, hence by \eqref{eq:forward-harnack},
$$
\sup_{B_\varrho} v(\cdot, t) \leq 2^{7+3\alpha}A \gamma r\omega(r).
$$
Since $t \in [-r^2\omega(r)^{2-p},0)$ was arbitrary, we obtain \eqref{eq:v-claim}.

With \eqref{eq:v-claim} at hand, we can now prove that
\begin{equation}\label{eq:v-claim2}
v(x,t) \leq \gamma\lambda r\omega(r),\quad\forall (x,t) \in P_r^- (\omega(r)),
\end{equation}
for some large constant $\lambda\equiv \lambda(n,p,\nu,L,A,\alpha)$, $\lambda > 1$, where $\gamma > 1$ again is the Harnack constant. To prove \eqref{eq:v-claim2}, it suffices to consider points $(x_0,t_0) \in P_r^-(\omega(r))$ such that
\begin{equation}\label{eq:v-claim2-re}
v(x_0,t_0) > \lambda r\omega(r).
\end{equation}
In what follows, we redefine $\theta$ as
$$
\theta \equiv \left(\frac{c}{v(x_0,t_0)}\right)^{p-2},
$$
where $c>1$ again is the waiting time constant in the intrinsic Harnack inequality. Let us choose $\lambda > c$ in \eqref{eq:v-claim2-re}. Then as $(x_0,t_0) \in P_r^-(\omega(r))$, we have
\begin{equation}\label{eq:v-claim2-re2}
t_0 + \theta |x_0|^p < t_0 + \left(\frac{c}{\lambda r\omega(r)}\right)^{p-2} |x_0|^p < t_0 \left[ 1 - \left(\frac{c}{\lambda}\right)^{p-2}\right] < 0.
\end{equation}
Next, we select $\lambda$ even larger if necessary such that
\begin{equation}\label{eq:v-claim2-re3}
t_0 - \theta (4|x_0|)^p > - \left[ 1 + 4^p \left( \frac{c}{\lambda}\right)^{p-2}\right] r^2 \omega(r)^{2-p} > -(8r)^2 \omega(8r)^{2-p}.
\end{equation}
Note that the first inequality in \eqref{eq:v-claim2-re3} is solely due to \eqref{eq:v-claim2-re} and  $(x_0,t_0) \in P_r^-(\omega(r))\subset Q_r^-(\omega(r))$. In view of \eqref{eq:v-claim2-re2} and \eqref{eq:v-claim2-re3}, we can employ the forward intrinsic Harnack inequality, which along with \eqref{eq:v-claim} yields that
\begin{equation}\label{eq:backward-harnack}
v(x_0,t_0) \leq \gamma \inf_{B_{|x_0|}} v(\cdot, t_0 + \theta_0|x_0|^p) \leq 2^{7+3\alpha} A\gamma^2 r\omega(r).
\end{equation}
Thus, \eqref{eq:v-claim2-re} implies \eqref{eq:backward-harnack}. Then by dichotomy, we arrive at \eqref{eq:v-claim2} with $\lambda \equiv 2^{7+3\alpha} A\gamma$. Finally, our assertion in the statement of the lemma follows by combining \eqref{eq:v-claim2} with \eqref{eq:v-u-deg-re}.
\end{proof}

Next, we establish the extension of the result in Lemma \ref{lem:C1a-deg-para} to entire cylinders under the additional assumption stated in Proposition \ref{prop:C1a-relower}.

\begin{lemma}\label{lem:C1a-deg-re}
Under the assumptions stated in Proposition \ref{prop:C1a-relower} (hence in particular we assume that \eqref{eq:criterion} holds for some
$\theta>0$), there exists $\lambda\equiv \lambda(n,p,\nu,L,A,\alpha,\theta)$ such that if $|D\ell| + s <\omega( r)<1$, then
$$
|u| \leq \lambda r \omega(r) \quad\text{in }Q_r(\omega(r)).
$$
\end{lemma}

\begin{proof} To prove the lemma we argue by way of contradiction. Indeed, let $H_j$, $s_j$, $\omega_j$, $u_j$, $\psi_j$, $\ell_j$ be as in the setting of this lemma, assume that $|D\ell_j| + s_j < \omega_j(r_j)$, that
\begin{equation}\label{eq:uj-C1a-re-asmp}
|u_j| < \lambda_j \varrho\omega(\varrho) \quad\text{in }Q_\varrho(\omega_j(\varrho)),
\end{equation}
for all $\varrho \in (r_j,1)$, for some $\lambda_j\to\infty$ and $r_j\to 0$, but that
\begin{equation}\label{eq:uj-C1a-re-false}
|u_j(\bar x_j, \bar t_j)| = \lambda_j r_j\omega_j(r_j),
\end{equation}
for some $(\bar x_j, \bar t_j)\in Q_{r_j}(\omega_j(r_j))$. That is, $\lambda_j$ is the value of the numerator in \eqref{eq:criterion} at $r = r_j$,  and hence
\begin{equation}\label{eq:uj-C1a-re-false-re}
|u_j(x_j,t_j)| \geq \theta^{-1}\lambda_j r_j \omega_j(r_j).
\end{equation}
for some $(x_j,t_j) \in Q_{r_j}(\lambda_j\omega_j(r_j))$.

As $r_j\to 0$, we can choose $\{\varrho_j\}$ to be a sequence satisfying
\begin{equation}\label{eq:vrj}
\lim_{j\to \infty}\varrho_j = 0,\quad \lim_{j\to\infty}\frac{\varrho_j}{r_j} = \infty,\quad\text{and}\quad \lim_{j\to\infty} \frac{1}{\lambda_j} \left(\frac{\varrho_j}{r_j}\right)^{1+\alpha} = 0.
\end{equation}
Let us now consider an auxiliary problem
\begin{equation}\label{eq:vj-pde-deg-re}
\begin{dcases}
  H_j   v_j = 0 &\text{in }Q_{\varrho_j}(\omega_j(\varrho_j)), \\
  v_j =   u_j + 2\varrho_j\omega_j(\varrho_j) &\text{on }\partial_{\rm P} Q_{\varrho_j}(\omega_j(\varrho_j)).
\end{dcases}
\end{equation}
Repeating, almost verbatim, the argument for the proof for \eqref{eq:v-u-deg} and \eqref{eq:v-u-deg-re}, we also deduce that
\begin{equation}\label{eq:vj-uj-deg-re}
\max\{0, u_j - 2\varrho_j\omega_j(\varrho_j)\} \leq v_j \leq u_j + 2\varrho_j \omega_j(\varrho_j) \quad\text{in }Q_{\varrho_j}(\omega_j(\varrho_j)).
\end{equation}
Now utilizing \eqref{eq:uj-C1a-re-asmp}, and considering large $j$, we have from \eqref{eq:vj-uj-deg-re} that
\begin{equation}\label{eq:vj-C1a-re-asmp-re}
0 \leq v_j \leq 2\lambda_j \varrho\omega_j(\varrho)\quad\text{in }Q_{\varrho_j}(\omega_j(\varrho)),
\end{equation}
for every $\varrho \in [r_j,\varrho_j)$, but by \eqref{eq:uj-C1a-re-false-re} and \eqref{eq:vrj},
\begin{equation}\label{eq:vj-C1a-re-false-re}
v_j(x_j,t_j) \geq {\theta^{-1}}\lambda_jr_j\omega_j(r_j) - 2\varrho_j\omega_j(\varrho_j) \geq \frac{1}{2\theta}\lambda_jr_j\omega_j(r_j),
\end{equation}
for all large $j$. On the other hand, by $u_j(0,0) = \psi_j(0,0) = 0$ (from assumption) and \eqref{eq:vj-uj-deg-re},
\begin{equation}\label{eq:vj-0-re}
v_j(0,0) \leq 2\varrho_j\omega_j(\varrho_j).
\end{equation}

We now introduce the change of variables
$$
(y,\tau) \equiv ({r_j}^{-1}x, r_j^{-2}{\theta_j^{p-2}t}),\quad \tilde v_j(y,\tau) \equiv \frac{v_j(x,t)}{r_j\theta_j},\quad\text{with}\quad \theta_j\equiv \lambda_j\omega_j(r_j),
$$
which maps $Q_{r_j}(\lambda_j\omega_j(r_j))$ onto $Q_1$, thus $Q_{r_j}(\omega_j(r_j))$ to $Q_1(\lambda_j^{-1})$. Since $\lambda_j\to \infty$, $Q_1(\lambda_j^{-1})$ approximates the infinite strip $B_1\times (-\infty,\infty)$. To simplify the notation, we let
$$
k_j \equiv \frac{\varrho_j}{r_j},\quad\text{and}\quad \tilde\omega_j (k) \equiv \frac{\omega_j(kr_j)}{\omega_j(r_j)}.
$$
Since $\omega_j$ verifies \eqref{eq:omega prop} with $A$ and $\alpha$, so does $\tilde\omega_j$.

We observe that
$$
\tilde H_j \tilde v_j = 0\quad\text{in }Q_{k_j}(\lambda_j^{-1}\tilde\omega_j(k_j)),
$$
where $\tilde H_j$ is a degenerate parabolic operator verifying \eqref{asp} with structure constants $\nu$, $L$ and the inhomogeneity constant $s \equiv s_j/\theta_j$. We remark that as $s_j \leq \omega_j(r_j)$ (from the setting), we have $s_j/\theta_j \to 0$ as $j\to\infty$.

By the compactness argument already detailed in the proof of Lemma \ref{lem:C1a-deg-int}, we can now deduce from \eqref{eq:vj-C1a-re-asmp-re} that $\tilde v_j \to \tilde v$ and  $D \tilde v_j \to D\tilde v$ locally uniformly in $\R^n\times \R$ along a subsequence (since any compact subset of $\R^n\times\R$ is contained in $Q_{k_j}(\lambda_j^{-1}\tilde\omega_j(k_j))$ for every large $j$), for some {\it nonnegative} weak solution $\tilde v$ to
$$
\tilde H \tilde v = 0\quad\text{in }\R^n\times\R.
$$
Writing by $\tilde a$ the vector-field associated with $\tilde H$, one can observe by passing to the limit in \eqref{asp} that $\langle \tilde a(z),z\rangle \geq\nu|z|^p$, $|\tilde a(z)| \leq L|z|^{p-1}$ and that $\langle \tilde a(z) - \tilde a(\xi), z- \xi\rangle > 0$ for every $z,\xi\in\R^n$, $z\neq\xi$.

Now we are in a position to derive a contradiction. Rephrasing \eqref{eq:vj-C1a-re-asmp-re} at $\varrho = r_j$ in terms of $\tilde v_j$ yields $0 \leq \tilde v_j \leq 2$ in $B_1\times (-\lambda_j^{-1},\lambda_j)$. By the locally uniform convergence, we obtain
$$
0\leq \tilde v \leq 2 \quad\text{in }B_1\times (-\infty,\infty).
$$
Using this it follows from a Liouville theorem for the degenerate parabolic equations, see Ch.\ 5, Proposition 5.3 and Remark 5.2 in \cite{DGV3}, that $\tilde v$ is constant in $\R^n\times\R$. However, we know from \eqref{eq:vj-0-re} and \eqref{eq:vrj} that $\tilde v_j (0,0) \to 0$, so the Liouville theorem ensures
$$
\tilde v\equiv 0\quad\text{in }\R^n\times\R.
$$
On the contrary, by \eqref{eq:vj-C1a-re-false-re} and recalling that $(x_j,t_j) \in Q_{r_j}(\lambda_j\omega_j(r_j))$, we have $\tilde v_j(y_j,\tau_j) \geq 1/(2\theta)$ for some $(y_j,t_j) \in Q_1$, so the uniform convergence yields $\tilde v \geq 1/(2\theta) > 0$, and this is a contradiction.

In conclusion, we observe that if $H$, $s$, $\omega$, $u$, $\psi$, $\ell$ are as in the setting of the lemma, then there exist  large $\lambda > 1$ and small $\e \in (0,1)$, such that if $|D\ell| + s < \e$, $|D\ell| + s < \omega(r) < 1$, and if
$$
|u| < \lambda \varrho \omega(\varrho) \quad\text{in }Q_\varrho(\omega(\varrho)),
$$
for every $\varrho \in (r,1)$, then the strict inequality continues to hold at $\varrho = r$. Thus, by the (standard) continuity argument previously outlined, the strict inequality holds whenever $|D\ell| + s < \varrho < 1$. This concludes the proof of the desired estimate.
\end{proof}

Next, we establish the extension of the result in Lemma \ref{lem:C1a-deg-para} to entire cylinders, without the additional assumption stated in Proposition \ref{prop:C1a-relower}, but only for $\alpha \in (0,\alpha_h)$, where $\alpha_h$ is as in Theorem \ref{thm:grad -}.

\begin{lemma}\label{lem:C1a-deg}
Under the assumptions stated in Proposition \ref{prop:C1alower}, there exists  $\lambda \equiv \lambda(n,p,\nu,L,A,\alpha)$ such that if  $|D\ell| + s < \omega(r) < 1$, then
$$
|u| \leq \lambda r \omega(r)\quad\text{in }Q_r(\lambda \omega(r)).
$$
\end{lemma}

\begin{proof}
The proof here is analogous with that of Lemma \ref{lem:C1a-deg-int}. In fact, most of the argument there can be repeated with $r_j^\alpha$ replaced with $\omega_j(r_j)$. The key difference is that we cannot play directly with $u_j$ as it is not a solution to the standard problem. For this reason, we shall construct an auxiliary function as in the proof of Lemma \ref{lem:C1a-deg-para}.

Let $H_j$, $s_j$, $u_j$, $\psi_j$, $\ell_j$, $\omega_j$ be as in the setting of this lemma such that $|D\ell_j| + s_j < \omega_j(r_j)$, and such that
\begin{equation}\label{eq:C1a-deg-asmp}
|u_j| < \lambda_j \varrho\omega(\varrho)\quad\text{in }Q_\varrho(\lambda_j\omega(\varrho)),
\end{equation}
for all $\varrho\in(r_j,1)$, for some $\lambda_j\to \infty$ and $r_j\to 0$, but
\begin{equation}\label{eq:C1a-deg-false}
|u_j(x_j,t_j)| = \lambda_jr_j\omega_j(r_j),
\end{equation}
for some $(x_j,t_j) \in Q_{r_j}(\lambda_j\omega_j(r_j))$.

As $r_j\to 0$, we can choose $\{\varrho_j\}$ as in \eqref{eq:vrj}. Consider the following auxiliary problem
\begin{equation}\label{eq:vj-pde-deg}
\begin{dcases}
  H_j   v_j = 0 &\text{in }Q_{\varrho_j}(\lambda_j\omega_j(\varrho_j)), \\
  v_j =   u_j + 2\varrho_j\omega_j(\varrho_j) &\text{on }\partial_{\rm P} Q_{\varrho_j}(\lambda_j\omega_j(\varrho_j)).
\end{dcases}
\end{equation}
Repeating, almost verbatim, the argument for the proof for \eqref{eq:v-u-deg} and \eqref{eq:v-u-deg-re} we also deduce that
\begin{equation}\label{eq:vj-uj-deg}
\max\{0, u_j - 2\varrho_j\omega_j(\varrho_j)\} \leq v_j \leq u_j + 2\varrho_j \omega_j(\varrho_j) \quad\text{in }Q_{\varrho_j}(\lambda_j\omega_j(\varrho_j)).
\end{equation}
Now utilizing \eqref{eq:C1a-deg-asmp}, we have from \eqref{eq:vj-uj-deg} that
\begin{equation}\label{eq:C1a-deg-asmp-re}
0 \leq v_j \leq 2\lambda_j \varrho\omega_j(\varrho)\quad\text{in }Q_\varrho(\lambda_j\omega_j(\varrho)),
\end{equation}
for every $\varrho \in (r_j,\varrho_j)$, but by \eqref{eq:C1a-deg-false} and \eqref{eq:vrj},
\begin{equation}\label{eq:C1a-deg-false-re}
v_j(x_j,t_j) \geq \lambda_jr_j\omega_j(r_j) - 2\varrho_j\omega_j(\varrho_j) \geq \frac{1}{2}\lambda_jr_j\omega_j(r_j),
\end{equation}
for all large $j$. On the other hand, by $u_j(0,0) = \psi_j(0,0) = 0$ (from the setting) and \eqref{eq:vj-uj-deg},
\begin{equation}\label{eq:vj-0}
v_j(0,0) \leq 2\varrho_j\omega_j(\varrho_j).
\end{equation}
By Theorem \ref{thm:grad -}, \eqref{eq:C1a-deg-asmp-re} (at $\varrho = \varrho_j$) implies
\begin{equation}\label{eq:Dvj}
|Dv_j(\cdot,t) - Dv_j(0,t)| \leq c \lambda_j\omega_j(\varrho_j)\left(\frac{r_j}{\varrho_j}\right)^{\alpha_h}\quad\text{in }B_{r_j},
\end{equation}
whenever $|t| \leq r_j^2 (\lambda_j\omega_j(\varrho_j))^{2-p}$.

With the last four displays at hand, we can now repeat the proof of Lemma \ref{lem:C1a-deg-int}. The only notable difference here is that we do not know {\it a priori} that $|D\tilde v_j(0,0)| \to 0$, where $\tilde v_j$ is the scaled version of $v_j$ as in \eqref{eq:tuj} with $\theta_j$ now given by $\lambda_j\omega_j(r_j)$. However, as $v_j$ is nonnegative and satisfies \eqref{eq:vj-0}, $\tilde v_j$ is also nonnegative and, due to \eqref{eq:vrj}, $\tilde v_j(0,0) \to 0$. Thus, the limit of $\{\tilde v_j\}$, if we call it $\tilde v$, satisfies $|D\tilde v(0,0)| = 0$, as $\tilde v \geq 0$ and $\tilde v(0,0) = 0$. We should also remark that as $\omega_j$ is assumed to verify \eqref{eq:omega prop} with $\alpha < \alpha_h$, \eqref{eq:Dvj} together with \eqref{eq:vrj} shows that $|D\tilde v_j(y,t) - D\tilde v_j(0,\tau)| \to 0$ for every $(y,\tau) \in Q_1$. Therefore, the same contradiction arises as in the proof of the aforementioned lemma. We omit further details.
\end{proof}


\subsection{The nondegenerate case}\label{sec:nondeg} We here study  the nondegenerate case, i.e., $|D\ell| + s \approx 1$. It is worth noting that in this case, the intrinsic $C^{1,\alpha}$ estimate holds for any $\alpha \leq {2}/{(p-2)}$. This will be proved by a linearization of the problem, and by making use of the strong minimum principle, the latter now being available due to the nondegeneracy of the obstacle.

\begin{lemma}\label{lem:C1a-ndeg}
Let $\alpha \in (0,2/(p-2)]$ in (place of $(0,\alpha_h)$) the assumptions in  Proposition \ref{prop:C1alower}, and assume in addition that $\mu\equiv |D\ell| + s > 0$ and that there exists
$\varrho\in (0,1)$,  $\omega(\varrho)\leq\mu$, such that
\begin{equation}\label{eq:C1a-ndeg-asmpprev}
|u-\ell| \leq \lambda\mu \varrho \quad\text{in }Q_\varrho(\lambda\mu),
\end{equation}
for some $\lambda > 1$. Then there exist $\bar\lambda\equiv \bar\lambda(n,p,\nu,L,A,\alpha,\lambda)$ and $c\equiv c(n,p,\nu,L)$, such that
\begin{equation}\label{eq:C1a-ndeg}
 |u-\ell| \leq \bar\lambda \mu r\frac{\omega(r)}{\omega(\varrho)} \quad\text{in }Q_r^-(c\lambda\mu),
\end{equation}
for every $r\in(0,\varrho)$.
\end{lemma}

\begin{remark}\label{rem:C1a-ndeg}
Our estimate in \eqref{eq:C1a-ndeg} depends on the two-sided bounds in \eqref{eq:C1a-ndeg-asmpprev} at the scale where the nondegeneracy is detected by the modulus of continuity of the obstacle. Such dependence cannot be removed, even at nondegenerate points, due to the presence of the self-similar solutions in \eqref{eq:barenblatt2}. This is also a key difference from the obstacle problems for uniformly parabolic equations, where the lower bound given by the obstacle implies the upper bound by the strong minimum principle.
\end{remark}

\begin{proof}[Proof of Lemma \ref{lem:C1a-ndeg}]
By \eqref{eq:p-C1a}, and that $\psi(0,0) = \ell(0) = 0$, $|D\ell| \leq \mu$, we have
$$
|\psi| \leq 2\mu r \quad\text{in }Q_r (\mu),
$$
for all $r\in(0,\varrho)$. Due to \eqref{eq:C1a-ndeg-asmpprev} and that $u(0,0) =\psi(0,0) = 0$, we see that
\begin{equation}\label{eq:C1a-ndeg-asmp}
|u| \leq \lambda\mu \varrho \quad\text{in }Q_\varrho(\lambda\mu).
\end{equation}
Furthermore, using also \cite[Lemma 4.4]{KMNobs} gives a large constant $c\equiv c(n,p,\nu,L)$ such that
\begin{equation}\label{eq:u-C01}
|u| \leq  c\lambda \mu r \quad\text{in }Q_r (c\lambda\mu),
\end{equation}
for every $r\in(0,\varrho)$. We shall use \eqref{eq:u-C01} later.

Let $r \in (0,\varrho/4)$ be given, and suppose that
\begin{equation}\label{eq:u-l-re}
|u - \ell | \leq 4\bar\lambda \mu r\frac{\omega(4r)}{\omega(\varrho)}\quad\text{in }Q_{4r}(c\lambda\mu),
\end{equation}
for some large constant $\bar\lambda\in(0,\frac{1}{2})$ to be determined later. Consider the auxiliary function problem
\begin{equation}\label{eq:v-pde-ndeg}
\begin{cases}
H v =  0 & \text{in }Q_{4r}(c\lambda\mu),\\
v = u + 4r\omega(4r) &\text{on }\partial_{\rm P}Q_{4r}(c\lambda\mu) .
\end{cases}
\end{equation}
Arguing as in Lemma \ref{lem:C1a-deg}, we can deduce from Lemma \ref{ellittico}, \eqref{eq:u-l-re} and $\omega(\varrho) < \mu$, that
\begin{equation}\label{eq:v-l}
0\leq v - \ell \leq 8\bar\lambda \mu r\frac{\omega(4r)}{\omega(\varrho)}\quad\text{in }Q_{4r}(c\lambda\mu),
\end{equation}
and that
\begin{equation}\label{eq:v-u}
|v - u | \leq 4r\omega(4r)\quad\text{in }Q_{4r}(c\lambda\mu).
\end{equation}

By \eqref{eq:u-C01}, \eqref{eq:v-u} and that $4r < \varrho < \omega^{-1}(\mu)$, we also have
$$
|v| \leq 4c\lambda \mu r\quad\text{in }Q_{4r}(c\lambda\mu).
$$
Hence, Theorem \ref{thm:grad -} implies that
\begin{equation}\label{eq:Dv-sup}
|Dv| \leq 4c_0c\lambda\mu\quad\text{in }Q_{2r}(c\lambda\mu),
\end{equation}
for some $c_0\equiv c_0(n,p,\nu,L)$.

We define an auxiliary function
\begin{equation}\label{eq:w}
w (x,t) \equiv  \frac{(v - \ell)(rx, r^2 (c\lambda\mu)^{2-p} t)}{8\bar\lambda \mu r\omega(4r)\omega(\varrho)^{-1}},
\end{equation}
which by \eqref{eq:v-l} satisfies
\begin{equation}\label{eq:w-re}
0 \leq  w \leq 1 \quad\text{in }Q_4.
\end{equation}
Making use of \eqref{eq:v-u} and that $u(0,0) = \psi(0,0) = \ell(0)$, we deduce that
\begin{equation}\label{eq:w-0}
w(0,0) \leq \frac{\omega(\varrho)}{2\mu \bar\lambda}.
\end{equation}
Moreover, as $Hv = H\ell = 0$ in $Q_{4r}(\lambda\mu)$, $w$ is a nonnegative solution to
\begin{equation}\label{eq:w-pde}
\ddiv (bDw)  - w_t = 0 \quad \text{in } Q_2,
\end{equation}
where $b(x,t)$, for each $(x,t)\in Q_4$, is a $(n\times n)$-dimensional matrix coefficient given by
\begin{equation}\label{eq:Br}
b(x,t)\equiv  \frac{1}{(c\lambda\mu)^{p-2}}\int_0^1 \partial_z a (\sigma D v(rx,r^2(c\lambda\mu)^{2-p}t)+ (1-\sigma)D\ell)\d\sigma.
\end{equation}
Utilizing \eqref{eq:Dv-sup} and that $\mu \equiv |D\ell| + s$, we deduce that
\begin{equation}\label{eq:b-ellip-up}
\sup_{Q_2}\langle b\xi,\xi\rangle \leq 2^{\frac{p-2}{2}}L|\xi|^2\left[ \frac{1}{(c\lambda\mu)^2} \sup_{Q_{2r}(c\lambda \mu)} |Dv|^2 +1\right]^{\frac{p-2}{2}} \leq c_3(c_2c_1)^{p-2} L|\xi|^2,
\end{equation}
for every $\xi\in\R^n$, for some $c_3\equiv c_3(p)$. The lower bound for the ellipticity of $b$ is more straightforward and we deduce
\begin{equation}\label{eq:b-ellip-low}
\begin{aligned}
\inf_{Q_2} \langle b\xi,\xi \rangle & \geq \frac{\nu}{(\lambda\mu)^{p-2}}|\xi|^2 \inf_{Q_{4r}(c\bar\lambda\mu)} \int_0^1 ( |\sigma Dv + (1-\sigma) D\ell|^2 + s^2)^{\frac{p-2}{2}}\d\sigma \\
& \geq \frac{c_4 \nu}{(c\lambda)^{p-2}}|\xi|^2 \left[ \frac{1}{\mu^2} \inf_{Q_{2r}(c\lambda\mu)} |Dv|^2 +1 \right]^{\frac{p-2}{2}}  \geq \frac{c_4  \nu}{(c\lambda)^{p-2}}|\xi|^2,
\end{aligned}
\end{equation}
for every $\xi\in\R^n$, for some constant $c_4 \equiv c_4(p)$, $0< c_4 < 1$.

Now given any $\e> 0$, we can employ Lemma \ref{lemma:Gaussian} to choose $\bar\lambda$ sufficiently large such that \eqref{eq:w-0} implies
\begin{equation}\label{ApplGaussian}
w \leq \e\quad\text{in }Q_1^-.
\end{equation}
In view of \eqref{eq:w-re}, \eqref{eq:b-ellip-up}, \eqref{eq:b-ellip-low} and $\omega(\varrho) < \mu$, the constant $\bar\lambda$ can be determined {\it a priori} by $n$, $p$, $\nu$, $L$, $\lambda$ and $\e$ only. Now rephrasing the last display in terms of $v - \ell$, we obtain
\begin{equation}\label{eq:v-l-re}
0 \leq v - \ell \leq 8\e\bar\lambda \mu r \frac{\omega(4r)}{\omega(\varrho)} \leq \frac{1}{2}\bar\lambda \mu  r\frac{\omega(r)}{\omega(\varrho)}\quad\text{in }Q_r^-(c\lambda\mu),
\end{equation}
where the derivation of the last inequality involves \eqref{eq:omega prop} to bound $\omega(4r)\leq 4^\alpha \omega(r)$, and choosing $\e$ sufficiently small so that $8\cdot 4^\alpha \e \leq \frac{1}{2}$. Note that $\e \equiv \e (A,\alpha)$, and hence that $\bar\lambda \equiv \bar\lambda ( n,p,\nu,L,A,\alpha,\lambda)$. Finally, selecting $\bar\lambda$ even larger such that $\bar\lambda > 2\cdot 4^{1+\alpha}$, it follows from \eqref{eq:v-u},\eqref{eq:omega prop} and \eqref{eq:u-l-re} that
\begin{equation}\label{eq:u-l}
|u - \ell| \leq \bar\lambda \mu r\frac{\omega(r)}{\omega(\varrho)} \quad\text{in }Q_r^-(c\lambda\mu).
\end{equation}

So far we have proved that \eqref{eq:u-l-re} implies \eqref{eq:u-l} whenever $r \in (0,\varrho/4)$. However, by \eqref{eq:C1a-ndeg-asmp}, \eqref{eq:u-l-re} holds with $r =\varrho/4$, provided that $\bar\lambda > \lambda$. Thus, we can iterate the implication for every sequence of radii $r =4^{-k}\varrho$, $k=1,2,\dots$, and prove that
$$
|u - \ell| \leq 4^{-k}\bar\lambda\mu \varrho\frac{\omega(4^{-k}\varrho)}{\omega(\varrho)}\quad\text{in }Q_{4^{-k}\varrho}^-(c\lambda\mu),
$$
for every $k =1,2,\dots$. The proof for \eqref{eq:C1a-ndeg} is now done by replacing $\bar\lambda$ with say $4^{1+\alpha}A\bar\lambda$ to cover the case when $r \in (4^{-k} \varrho, 4^{-k+1}\varrho)$.
\end{proof}

Without any information on the upper bound, the best we obtain here is the linear development of the solution in the intrinsic backward $p$-paraboloids. It is unclear, however, whether such an estimate can be improved via linearization.

\begin{lemma}\label{lem:C1a-ndeg-para}
Let $\alpha \in (0,2/(p-2)]$ in (place of $(0,\alpha_h)$) the assumptions in  Proposition \ref{prop:C1alower}, and assume in addition that $\mu\equiv |D\ell| + s > 0$. Then there exists $\lambda\equiv \lambda(n,p,\nu,L)$ such that
$$
|u| \leq \lambda \mu r\quad\text{in }P_r^-(\mu),
$$
for every $r$ such that $0<\omega(r)<\mu$.
\end{lemma}

\begin{proof}
The proof is essentially the same as that of Lemma \ref{lem:C1a-deg-para}; we only need to replace $\omega(r)$ there with $\mu$. Hence, we omit the details.
\end{proof}

\begin{remark} Note that it is only in the argument leading to \eqref{ApplGaussian}, as we rely on Lemma \ref{lemma:Gaussian} and the strong minimum principle, that our argument needs a restriction to  lower intrinsic cylinders.\end{remark}


\subsection{Proofs of Propositions \ref{prop:C1a-para}-\ref{prop:C1alower}}
The proofs of  Proposition \ref{prop:C1a-para}, Proposition \ref{prop:C1a-relower}, and Proposition \ref{prop:C1alower}, are essentially the same as that of Proposition \ref{prop:C1a-int}. We first apply Lemma \ref{lem:C1a-deg-para}, Lemma \ref{lem:C1a-deg-re}, and Lemma \ref{lem:C1a-deg}, respectively, in place of Lemma \ref{lem:C1a-deg-int}, to the scale of degeneracy, $\mu \equiv |D\ell| + s$. If $\mu = 0$, then we are done. If $\mu < 1$, then we rescale our problem so that the scaled problem is no longer degenerate. Then, we apply Lemma \ref{lem:C1a-ndeg-para}, Lemma \ref{lem:C1a-ndeg}, and Lemma \ref{lem:C1a-ndeg}, respectively, in place of Lemma \ref{lem:C1a-ndeg-int}, to finish the proofs. Note that if $\mu < 1$, then we apply Lemma \ref{lem:C1a-ndeg} in the proof of Proposition \ref{prop:C1a-relower} as well as in the proof of Proposition \ref{prop:C1alower}. We leave out the details.


\subsection{Proofs of Theorem \ref{thm:C1a}}\label{pproffa}

The theorem follows immediately from the following rescaling argument. We introduce the change of variables
$$
(y,\tau) \equiv ( \bar r^{-1}{(x-x_0)}, {\bar r^{-2}\bar\mu^{p-2}(t-t_0)}),
$$
which maps $Q_{\bar r}(\bar \mu)$ onto $Q_1$. We introduce the auxiliary functions
$$
\tilde u (y,\tau) \equiv \frac{u(x,t) - \psi(x_0,t_0)}{\bar\mu\bar r}\quad\text{and}\quad \tilde \psi(y,\tau) \equiv \frac{\psi(x,t) - \psi(x_0,t_0)}{\bar\mu\bar r}.
$$
Then $\tilde u (0,0) = \psi(0,0) = 0$. Moreover,
$$
\max\{ \tilde H \tilde u, \tilde\psi - \tilde u\} = 0\quad\text{in }Q_1,
$$
where $\tilde H\tilde u \equiv \ddiv \tilde a(D\tilde u) - \tilde u_\tau$, with $\tilde a (D\tilde u (y,\tau)) \equiv \bar\mu^{1-p} a(D u(x,t))$. Hence, $\tilde H$ verifies \eqref{asp} but with the inhomogeneity parameter $\tilde s \equiv s/\bar\mu$. It is now easy to verify that $\tilde H$, $\tilde s$, $\tilde u$ and $\tilde\psi$ verify the assumptions of Proposition \ref{prop:C1a-para}. The final estimate then follow immediately by rewriting the assertion of the proposition in terms of $u$.



\section{Regularity up to the initial layer/state}\label{sec:initial}

This section is devoted to the study of the regularity estimate at the initial (time) layer/state.

\begin{proposition}\label{prop:C1a-ini}
Let $H$ be as in \eqref{generalH}, with $(a,s)$ as in \eqref{asp}, let $\omega$ verify \eqref{eq:omega prop} with $\alpha\leq {2}/{(p-2)}$, and let $u$ be a continuous weak solution to
$$
\begin{cases}
\max\{ Hu , \psi - u\} = 0 & \text{in } Q_1^+, \\
u = g & \text{on } B_1\times \{0\},
\end{cases}
$$
such that $u(0,0) = g(0) = 0$ and $|u|\leq 1$ in $Q_1^+\cup\partial_{\rm P}Q_1^+$. Let $\ell$ and $\tilde\ell$ be time-independent affine functions, such that $\max\{|D\ell|, |D\tilde\ell|\} + s \leq 1$, and set $\vp(r) \equiv \max\{\omega(r), |D\ell| + s\}$ and $\tilde\vp(r) \equiv \max\{\omega(r), |D\tilde\ell|+s\}$. Assume that
\begin{equation}\label{eq:g-C1a}
 |g - \ell| \leq r\omega(r)\quad\text{in }B_r,\quad | \psi - \tilde\ell| \leq r\omega(r)\quad\text{in }Q_r^+(\tilde\vp(r)),
\end{equation}
for every $r\in(0,1)$. Then for some $\lambda\equiv \lambda(n,p,\nu,L,A,\alpha)$,
\begin{equation}\label{eq:C1a-ini}
 |u - \ell | \leq \lambda r\omega(r)\quad\text{in }Q_r^+(\lambda\vp(r)),
\end{equation}
for every $r\in(0,1)$.
\end{proposition}

\begin{remark}\label{rem:C1a-ini}
Note that this proposition allows for the case $g(0 ) > \psi(0,0)$, that is, the reference point does not need to be a contact point.
\end{remark}

The  proof of Proposition \ref{prop:C1a-ini} will be presented in subsections and we again divide the argument into the degenerate case ($|D\ell| + s \approx 0$), and the nondegenerate case ($|D\ell| + s \approx 1$).


\subsection{The degenerate case} We here study the degenerate case, $|D\ell| + s \approx 0$. The idea is to resort to the uniqueness result for degenerate parabolic Cauchy problem in \cite{KuPa}. Our proof leaves $|D\tilde\ell| + s$ free, hence it also works when $|D\tilde\ell| + s \approx 1$.

\begin{lemma}\label{lem:C1a-ini-deg} Assume, in addition to the assumptions in Proposition \ref{prop:C1a-ini}, that $|D\ell| +s < 1$. Then there exists $\lambda\equiv \lambda(n,p,\nu,L,A,\alpha)$, $\lambda > 1$, such that if $|D\ell| + s < \omega(r) < 1$, then
$$
|u| \leq \lambda r\omega(r)\quad\text{in }Q_r^+(\lambda\omega(r)).
$$
\end{lemma}

\begin{proof}
For each $j=1,2,\dots$, let $H_j$, $s_j$, $\omega_j$, $\ell_j$, $\tilde\ell_j$, $g_j$, $\psi_j$ and $u_j$ be as in the setting of the lemma, $|D\ell_j| + s_j \leq \omega_j(r_j)$, and
\begin{equation}\label{eq:C1a-ini-deg-asmp}
|u_j| < \lambda_j \varrho \omega_j(\varrho)\quad\text{in }Q_\varrho^+(\lambda_j\omega_j(\varrho)),
\end{equation}
for all $\varrho \in (r_j,1)$, for some $\lambda_j\to\infty$ and $r_j\to 0$, but
\begin{equation}\label{eq:C1a-ini-deg-false}
|u_j(x_j,t_j)| = \lambda_j r_j\omega_j(r_j),
\end{equation}
for some $(x_j,t_j) \in Q_{r_j}^+(\lambda_j\omega_j(r_j))\cup\partial_{\rm P}Q_{r_j}^+(\lambda_j\omega_j(r_j))$.

Since $r_j\to 0$,   we can choose a sequence $\{\varrho_j\}$ in $(0,\omega_j^{-1}(\mu))$  such that $\{\varrho_j\}$ satisfies \eqref{eq:vrj}. Let $\eta_j$ be a smooth cutoff function (in space only) verifying $0\leq \eta_j \leq 1$ in $B_{\varrho_j}$, $\eta_j = 1$ in $B_{\varrho_j/2}$, and $\supp(\eta_j)\subset B_{3\varrho_j/4}$. Using $\eta_j$, we introduce the auxiliary problem
\begin{equation}\label{eq:vj-pde-ini-deg}
\begin{dcases}
  H_j   v_j = 0 &\text{in }Q_{\varrho_j}^+(\lambda_j\omega_j(\varrho_j)), \\
  v_j =   u_j + 5\varrho_j\omega_j(\varrho_j) &\text{on }\partial_{\rm P} Q_{\varrho_j}^+(\lambda_j\omega_j(\varrho_j))\setminus(B_{\varrho_j}\times\{0\}), \\
  v_j = (1-\eta_j)   g_j +  5\varrho_j\omega_j(\varrho_j) & \text{on }B_{\varrho_j}\times\{0\}.
\end{dcases}
\end{equation}
Since $\supp\eta_j\subset B_{3\varrho_j/4}$, $v_j$ is continuous on $\partial_{\rm P} Q_{\varrho_j}^+(\lambda_j\omega_j(\varrho_j))$. We claim that
\begin{equation}\label{eq:vj-uj-ini-deg}
0\leq   v_j -   u_j \leq 5\varrho_j\omega_j(\varrho_j)\quad\text{in }Q_{\varrho_j}^+(\lambda_j\omega_j(\varrho_j)).
\end{equation}
Since
$$
  v_j \leq   u_j + 5\varrho_j\omega_j(\varrho_j) \quad\text{on }\partial_{\rm P}Q_{\varrho_j}^+(\lambda_j\omega_j(\varrho_j)),
$$
and
$$
  H_j   u_j \leq 0\quad\text{in }Q_{\varrho_j}^+(\lambda_j\omega_j(\varrho_j)),
$$
the upper bound follows immediately from the standard comparison principle. The lower bound needs extra care, and we will again, as in the proof of Lemma \ref{lem:C1a-deg}, invoke the ``elliptic'' comparison principle, Lemma \ref{ellittico}.

First, observe from \eqref{eq:g-C1a}, that $g_j(0) = 0$ and $\omega_j(\varrho_j) > \omega_j(r_j) \geq |D\ell_j|$, that
$$
  v_j \geq g_j + 4\varrho_j \omega_j(\varrho_j) \geq \ell_j + 3\varrho_j\omega_j(\varrho_j)  \quad\text{on }B_{\varrho_j}\times\{0\}.
$$
However, as $g_j\geq \psi_j$ on $B_{\varrho_j}\times \{0\}$, it follows from \eqref{eq:g-C1a} that
$$
\ell_j - \tilde\ell_j \geq (\ell_j - g_j) + (\psi_j -\tilde\ell_j) \geq - 2\varrho_j\omega_j(\varrho_j)\quad\text{on }B_{\varrho_j}\times\{0\}.
$$
Combining the last two displays, we obtain
$$
v_j \geq \tilde\ell_j + \varrho_j\omega_j(\varrho_j)\quad\text{on }B_{\varrho_j}\times\{0\}.
$$
However, by $  u_j\geq   \psi_j$ in $\partial_{\rm P}Q_{\varrho_j}^+(\lambda_j\omega_j(\varrho_j))$ and \eqref{eq:g-C1a} again, we also see that
$$
  v_j \geq  \psi_j + 4\varrho_j\omega_j(\varrho_j) \geq  \tilde\ell_j + 3\varrho_j\omega_j(\varrho_j)\quad\text{on }\partial_{\rm P} Q_{\varrho_j}^+(\lambda_j\omega_j(\varrho_j))\setminus (B_{\varrho_j}\times\{ 0\}).
$$
By the comparison principle and the evident fact that $  H_j \ell_j \equiv 0$, we deduce that
\begin{equation}\label{eq:vj-lj'}
v_j \geq \tilde\ell_j + \varrho_j\omega_j(\varrho_j)\quad\text{on }Q_{\varrho_j}^+(\lambda_j\omega_j(\varrho_j)).
\end{equation}

With \eqref{eq:vj-lj'} at hand, we may use \eqref{eq:g-C1a} again to derive that
$$
  v_j \geq   \psi_j\quad\text{in }Q_{\varrho_j}^+(\lambda_j\omega_j(\varrho_j)),
$$
in particular, $  v_j \geq   u_j$ on $Q_{\varrho_j}^+(\lambda_j\omega_j(\varrho_j))\cap\{  u_j =  \psi_j\}$. So far we have observed that
$$
  v_j \geq   u_j\quad\text{on }\partial (\{   u_j >  \psi_j\}\cap Q_{\varrho_j}^+(\lambda_j\omega_j(\varrho_j))),
$$
hence by Lemma \ref{ellittico}, and that $  H_j   v_j =   H_j   u_j = 0$ in the open set $\{  u_j >   \psi_j\}\cap Q_{\varrho_j}^+(\lambda_j\omega_j(\varrho_j))$, we obtain the lower bound in \eqref{eq:vj-uj-ini-deg}. Thus, both inequalities in \eqref{eq:vj-uj-ini-deg} are proved.

We introduce
\begin{equation}\label{eq:kj-ej-ini}
k_j \equiv  \frac{\varrho_j}{r_j}, \quad  \e_j\equiv  \frac{k_j\omega_j(\varrho_j)}{\lambda_j\omega_j(r_j)},\quad \tilde\omega_j(k) \equiv  \frac{\omega_j(kr_j)}{\omega_j(r_j)},
\end{equation}
and the scaled version of $v_j$,
$$
\tilde v_j (y,\tau) \equiv  \frac{v_j(r_jy , r_j^2 \theta_j^{2-p}\tau)}{r_j\theta_j},\quad\text{with}\quad \theta_j \equiv  \lambda_j \omega_j(r_j).
$$
Recalling that $\eta_j\equiv 1$ in $B_{\varrho_j/2}$, it follows from \eqref{eq:vj-pde-ini-deg} that $\tilde v_j$ is a weak solution to
\begin{equation}\label{eq:tvj-pde-ini-deg}
\begin{dcases}
\tilde H_j \tilde v_j = 0 &\text{in }Q_{k_j}^+(\tilde\omega_j(k_j)), \\
\tilde v_j = 4\e_j &\text{on }B_{k_j/2} \times \{0\},
\end{dcases}
\end{equation}
where $\tilde H_j$ is the degenerate parabolic operator verifying \eqref{asp} with the inhomogeneity constant $s\equiv s_j/\theta_j$. We remark that $s_j/\theta_j \to 0$, since $s_j \leq \omega_j(r_j)$ and $\lambda_j\to \infty$.

Utilizing \eqref{eq:C1a-ini-deg-asmp}, \eqref{eq:vrj} and \eqref{eq:kj-ej-ini}, we can deduce from  \eqref{eq:vj-uj-ini-deg} that
\begin{equation}\label{eq:tvj-osc-ini}
\begin{aligned}
\osc_{Q_k^+(\tilde\omega_j(k))} \tilde v_j \leq k\tilde\omega_j(k) + \frac{5k_j\tilde\omega_j(k_j)}{\lambda_j} \leq 2k\tilde\omega_j(k),
\end{aligned}
\end{equation}
for every $k\in (1,k_j)$ and all large $j$. Moreover, by \eqref{eq:C1a-ini-deg-false}, \eqref{eq:vrj} and \eqref{eq:kj-ej-ini},
\begin{equation}\label{eq:tvj-osc-ini-re}
|\tilde v_j (y_j,\tau_j)| \geq 1 - \frac{5k_j\tilde\omega_j( k_j)}{\lambda_j} \geq\frac{1}{2},
\end{equation}
for some $(y_j,\tau_j) \in Q_1^+\cup\partial_{\rm P}Q_1^+$.

As $\tilde\omega$ also verifies \eqref{eq:omega prop}, with $\alpha\leq{2}/{(p-2)}$, we have
$$
Q_1^+ \subset B_k\times (0,1) \subset Q_k^+(\tilde\omega_j(k)),\quad\forall  j,k = 1,2,\dots.
$$
By \eqref{eq:tvj-pde-ini-deg} (it is worth noting that $\tilde v_j$ is constant on the initial layer, $B_{k_j/2}\times\{0\}$ and the constant, $\e_j$, converges to zero as $j\to\infty$ by \eqref{eq:kj-ej-ini} and \eqref{eq:vrj}) and \eqref{eq:tvj-osc-ini}, it follows from \cite{DBb} and \cite{Li} that $\{\tilde v_j\}$ and $\{D\tilde v_j\}$ are uniformly H\"older continuous families in $B_k\times [0,1]$, for each $k=1,2,\dots$. Hence, by compactness, $\tilde v_j \to \tilde v$  and $D\tilde v_j \to D\tilde v$ uniformly in $B_k\times [0,1]$, along a subsequence, for some continuous function $\tilde v:\R^n\times[0,1]\to \R$ having continuous spatial derivatives.

Let $\tilde a_j$ denote the symbol associated with $\tilde H_j$. Recall that $\tilde a_j$ verifies \eqref{asp} with the inhomogeneity parameter $s = s_j/\theta_j$ and that $s_j/\theta_j \to 0$. Thus, extracting a further subsequence along which $D\tilde v_j\to D\tilde v$, we have $\tilde a_j \to \tilde a$ locally uniformly in $\R^n$, where $\tilde a$ is a Lipschitz vector-field on $\R^n$. Passing to the limit in \eqref{asp} with $s = s_j/\theta_j \to 0$, we observe that $\tilde a$ is homogeneous, strictly monotone and has $p$-growth, i.e., $\langle \tilde a(\xi),\xi\rangle\geq\nu|\xi|^p$, $|\tilde a(\xi)| \leq L|\xi|^{p-1}$, and $\langle a(z)- a(\xi), z-\xi\rangle > 0$ for all $z,\xi\in\R^n$ with $z\neq \xi$.

From the above observations, we have $\tilde a_j(D\tilde v_j) \to \tilde a(D\tilde v)$ locally uniformly in $\R^n\times(0,1)$. Thus, we may pass to the limit in the weak formulation of \eqref{eq:tvj-pde-ini-deg}, use $\e_j\to 0$ as well as the uniform convergence of $\tilde v_j\to \tilde v$ in every compact subset of $\R^n\times[0,1]$, to derive that
$$
\begin{dcases}
\tilde H\tilde v = 0 &\text{in }\R^n\times (0,1),\\
\tilde v = 0 & \text{on }\R^n\times\{0\},
\end{dcases}
$$
where $\tilde H$ is the degenerate parabolic operator associated with the symbol $\tilde a$. Since $\tilde a$ is homogeneous, strictly monotone and has $p$-growth, we can employ the uniqueness result for  degenerate parabolic Cauchy problems, \cite[Theorem 5.4]{KuPa}, to infer that $\tilde v \equiv 0$ in $\R^n\times[0,1)$. Nevertheless, letting $j\to \infty$ in \eqref{eq:tvj-osc-ini-re} and using that $\tilde v_j\to \tilde v$ uniformly in $\overline{Q_1^+} \equiv \overline{B_1}\times[0,1]$, we find that
$$
\tilde v(y_0,\tau_0) \geq \frac{1}{2},
$$
for some $(y_0,\tau_0)\in \overline{Q_1^+}\equiv \overline{B_1}\times[0,1]$. Thus, we have arrived at a contradiction.

What we have proved is that there is a large constant $\bar\lambda\equiv \bar\lambda(n,p,\nu,L,A,\alpha)$, and a small constant $\bar r\equiv \bar r(n,p,\nu,L,A,\alpha)$, such that  if for some $\lambda > \bar\lambda$,  and for some $r < \bar r$ satisfing $|D\ell| + s < \omega(r)$, we have
$$
|u| < \lambda \varrho \omega(\varrho)\quad\text{in }Q_\varrho^+(\lambda\omega(\varrho)),
$$
for all $\varrho\in(r,1)$, then the strict inequality continues to hold at $\varrho = r$. The desired estimate then follows immediately by the continuity argument used before.
\end{proof}


\subsection{The nondegenerate case}\label{sec:ini-gen} We here study the nondegenerate case, $ |D\ell| + s \approx 1$. Our argument is based on the uniqueness result for uniformly parabolic Cauchy problem with unbounded measurable coefficients. Such Cauchy problems are induced by linearizing our solution around approximating affine functions. We remark that also this proof leaves the (non-)degeneracy of the obstacle $\psi$ free, and in particular, it allows for $|D\tilde\ell| + s \approx 0$.

\begin{lemma}\label{lem:C1a-ini-ndeg}Assume, in addition to the assumptions in  Proposition \ref{prop:C1a-ini}, that $\mu\equiv |D\ell| + s > 0$ and that
\begin{equation}\label{eq:C1a-ini-ndeg-asmp}
|u| \leq \lambda \mu \varrho\quad\text{in }Q_\varrho^+(\lambda\mu),
\end{equation}
for some $\varrho \in (0,\omega^{-1}(\mu))$ and some $\lambda > 1$. Then there exist $\bar\lambda\equiv\bar\lambda(n,p,\nu,L,A,\alpha,\lambda)$ and $c\equiv c(n,p,\nu,L)$, $\bar\lambda >\lambda$ and $c > 1$, such that
\begin{equation}\label{eq:C1a-ini-ndeg}
|u - \ell| \leq \bar\lambda \mu r \frac{\omega(r)}{\omega(\varrho)}\quad\text{in }Q_r^+(c\lambda\mu),
\end{equation}
for every $r\in (0,\varrho)$.
\end{lemma}

\begin{proof}
The proof is similar to that of Lemma \ref{lem:C1a-ini-deg} in the sense that we will resort to the uniqueness for Cauchy problems. The only difference here is that due to the nondegeneracy of the initial data ($g_j$ below), the Cauchy problems are indeed uniformly parabolic.

Let $\lambda > 1$ be given, and let $H_j$, $s_j$, $\omega_j$, $u_j$, $g_j$, $\psi_j$, $\ell_j$ and $\tilde\ell_j$ be  as in the lemma, and such that
\begin{equation}\label{eq:C1a-ini-ndeg-asmp-re2}
|u_j| \leq \lambda \mu_j \bar\varrho_j \quad\text{in }Q_{\bar\varrho_j}^+(\lambda \mu_j),
\end{equation}
for some $\bar\varrho_j \in (0,\omega_j^{-1}(\mu_j))$ with $\mu_j \equiv |D\ell_j| + s_j > 0$. By \eqref{eq:g-C1a} and \eqref{eq:C1a-ini-ndeg-asmp}, \cite[Lemma 4.2]{KMNobs} yields that
\begin{equation}\label{eq:uj-osc-ini-re}
|u_j| \leq c\lambda \mu_j \varrho \quad\text{in }Q_\varrho^+(c\lambda\mu_j),
\end{equation}
for every $\varrho \in (0,\bar\varrho_j)$, for some $c\equiv c(n,p,\nu,L)$, $c> 1$.

We now assume, by way of contradiction, that
\begin{equation}\label{eq:C1a-ini-ndeg-asmp-re}
|u_j - \ell_j | < \lambda_j \mu_j \varrho \frac{\omega_j(\varrho)}{\omega_j(\bar\varrho_j)} \quad\text{in }Q_\varrho^+(\lambda\mu_j),
\end{equation}
for every $\varrho \in (r_j, \bar\varrho_j)$, for some $\lambda_j\to \infty$, $\lambda_j > \lambda$ and $r_j \to 0$, but that
\begin{equation}\label{eq:C1a-ini-ndeg-false}
|(u_j - \ell_j)(x_j,t_j) | \geq \lambda_j \mu_j r_j \frac{\omega_j(r_j)}{\omega_j(\bar\varrho_j)},
\end{equation}
for some $(x_j,t_j) \in Q_{r_j}^+(c\lambda\mu_j)$. Once we derive the contradiction, we can argue as at the end of the proof of Lemma \ref{lem:C1a-ini-deg} to deduce that \eqref{eq:C1a-ini-ndeg-asmp} implies \eqref{eq:C1a-ini-ndeg}.

Since $r_j\to 0$, we can choose a sequence $\{\varrho_j\}$ in $(0,\bar\varrho_j)$ satisfying \eqref{eq:vrj}. We let $\eta_j$ be as in the proof of Lemma \ref{lem:C1a-ini-deg}, and we consider the auxiliary problem
\begin{equation}\label{eq:vj-pde-ini}
\begin{dcases}
  H_j   v_j = 0 &\text{in }Q_{\varrho_j}^+(c\lambda\mu_j), \\
  v_j =   u_j + 5\varrho_j\omega_j(\varrho_j)&\text{on }\partial_{\rm P} Q_{\varrho_j}^+(c\lambda\mu_j)\setminus(B_{\varrho_j}\times\{0\}), \\
  v_j = (1-\eta_j)   g_j + \eta_j \ell_j + 5\varrho_j\omega_j(\varrho_j)& \text{on }B_{\varrho_j}\times\{0\}.
\end{dcases}
\end{equation}
Arguing as in the proof of Lemma \ref{lem:C1a-ini-deg}, we deduce that
\begin{equation}\label{eq:vj-uj-ini-ndeg}
0\leq   v_j -   u_j \leq 5\varrho_j\omega_j(\varrho_j)\quad\text{in }Q_{\varrho_j}^+(c\lambda\mu_j).
\end{equation}
On the one hand, together with \eqref{eq:C1a-ini-ndeg-asmp-re}, \eqref{eq:vj-uj-ini-ndeg} and \eqref{eq:omega prop}, \eqref{eq:vrj} implies that
\begin{equation}\label{eq:vj-lj-ini}
\begin{aligned}
|v_j - \ell_j| &\leq \lambda_j \varrho\omega_j(\varrho) + 5\varrho_j\omega_j(\varrho_j) \leq  2A\lambda_j \omega_j(r_j) \varrho \left(\frac{\varrho}{r_j}\right)^\alpha\quad\text{in }Q_\varrho^+(c\lambda\mu_j),
\end{aligned}
\end{equation}
for every $\varrho\in(r_j,\varrho_j)$ and for all large $j$. In addition, it follows from \eqref{eq:C1a-ini-ndeg-false} and \eqref{eq:vj-uj-ini-ndeg} that
\begin{equation}\label{eq:vj-lj-ini-re}
 |(v_j - \ell_j)(x_j,t_j)| \geq \lambda_j r_j \omega_j(r_j)  -  5\varrho_j\omega_j(\varrho_j) \geq \frac{1}{2}\lambda_j r_j \omega_j(r_j),
\end{equation}
again for all large $j$. On the other hand, \eqref{eq:uj-osc-ini-re}, along with $\omega_j(\varrho_j) < \omega_j(\bar\varrho_j) < \mu_j$ and that $\lambda,c > 1$,  implies that $|v_j| \leq 2c\lambda \mu_j\varrho_j$ in $Q_{\varrho_j}^+(c\lambda\mu_j)$, and hence Theorem \ref{thm:grad -} yields that
\begin{equation}\label{eq:Dvj-ini-ndeg}
|Dv_j| \leq c_0c\lambda\mu_j \quad\text{in }Q_{\varrho_j/4}^+(c\lambda\mu_j).
\end{equation}

Let $k_j$ and $\e_j$ be as in \eqref{eq:kj-ej-ini} and define
\begin{equation}\label{eq:kj-ej-iniuu}
w_j (y,\tau) \equiv  \frac{(v_j - \ell_j) (r_j y, r_j^2 (c\lambda\mu_j)^{2-p} \tau)}{\lambda_jr_j\omega_j(r_j)}.
\end{equation}
Since $H_j v_j = H_j \ell_j = 0$ in $Q_{\varrho_j}^+(c\lambda\mu_j)$ and $\eta_j \equiv 1$ on $B_{\varrho_j/2}$, $w_j$ is a continuous weak solution to
\begin{equation}\label{eq:wj-pde-ini}
\begin{dcases}
P_j w_j = 0 &\text{in }Q_{k_j}^+, \\
w_j = 5\e_j & \text{on }B_{k_j/2} \times \{0\},
\end{dcases}
\end{equation}
where $P_j$ is the linear parabolic operator,
$$
P_j w \equiv \ddiv (b_j D w) - w_\tau,
$$
associated with the matrix-valued map $b_j : Q_{k_j}^+ \to \R^{n\times n}$,
$$
b_j (y,\tau) \equiv \frac{1}{(c\lambda\mu_j)^{p-2}}\int_0^1 \partial_z a_j (\sigma Dv_j(r_jy,r_j^2 (c\lambda\mu_j)^{2-p} \tau) + (1-\sigma) D\ell_j)\,\d\sigma.
$$
By \eqref{asp}, \eqref{eq:Dvj-ini-ndeg} and that $\mu_j\equiv |D\ell_j| + s_j$, we readily deduce (as in the derivation of  \eqref{eq:b-ellip-up} and \eqref{eq:b-ellip-low}) that
\begin{equation}\label{eq:bj-ellip-ini}
\frac{c_1 \nu}{(c\lambda)^{p-2}}|\xi|^2 \leq  \langle b_j\xi,\xi\rangle  \leq c_2c_0^{p-2} L |\xi|^2 \quad\text{in }Q_{k_j/4}^+,
\end{equation}
for all $\xi\in\R^n$, for some $c_i \equiv c_i(p)$, $i\in\{1,2\}$.

We remark that \eqref{eq:vj-lj-ini} yields compactness of $\{w_j\}$. In fact, by \eqref{eq:vj-lj-ini}, $|w_j| \leq 2Ak^{1+\alpha}$ in $Q_k^+$ for every $k\in(1,k_j)$. Now by \eqref{eq:wj-pde-ini} and \eqref{eq:bj-ellip-ini}, it follows from \cite{Lieberman_book} that
\begin{equation}\label{eq:wj-w-ini}
\begin{dcases}
Dw_j \to Dw  \quad\text{weakly in }\mathrm{L}^2 (0,k^2; \mathrm{L}_{loc}^2(B_k)), \\
w_j \to w \quad\text{uniformly in }\overline{B_k} \times [0,k^2],
\end{dcases}
\end{equation}
along a subsequence as $j\to\infty$, for each $k=1,2,\dots$, for some $w\in \mathrm{L}^2(0,\infty; W_{loc}^{1,2}(\R^n)) \cap C(\R^n\times[0,\infty))$. By \eqref{eq:wj-w-ini} and \eqref{eq:bj-ellip-ini}, we can pass to the limit in the weak formulation (by the parabolic $G$-convergence, see \cite{Spa}) of \eqref{eq:wj-pde-ini} to derive that
$$
\begin{dcases}
P w = 0 &\text{in } \R^n\times(0,\infty), \\
w = 0 & \text{on } \R^n\times \{0 \},
\end{dcases}
$$
where $P$ is a linear parabolic operator associated with uniformly parabolic and bounded measurable coefficients in the whole space $\R^n\times(0,\infty)$. Since $|w_j|\leq 2Ak^{1+\alpha}$ in $Q_k^+$ for every $k\in(1,k_j)$ and all large $j$, the uniform convergence implies $|w| \leq 2Ak^{1+\alpha}$ in $Q_k^+$ for every $k=1,2,\dots$, i.e., $w$ has polynomial growth at infinity. As $w = 0$ on the initial layer $\tau = 0$, the uniqueness result for uniformly parabolic linear Cauchy problems, see \cite{AB}, ensures $w\equiv 0$ in the whole space $\tau \geq 0$. This however yields a contradiction. Indeed, rephrasing \eqref{eq:vj-lj-ini-re} in terms of $w_j$ yields
$$
|w_j(y_j,\tau_j)| \geq \frac{1}{2},
$$
for some $(y_j,\tau_j) \in Q_1^+$, and thus letting $j\to\infty$ and using the uniform convergence of $w_j\to w$ yields $w \not\equiv 0$ in $\overline{Q_1^+}$.
\end{proof}


\subsection{Proof of Proposition \ref{prop:C1a-ini}}\label{sec:pf-ini}

The idea of the proof here is the same as  that of Proposition \ref{prop:C1alower}. In fact, one may follow the lines of the proof, replacing the full space-time cylinders with forward space-time cylinders, and use Lemmas \ref{lem:C1a-ini-ndeg} and \ref{lem:C1a-ini-deg} in place of Lemmas \ref{lem:C1a-ndeg} and \ref{lem:C1a-deg} respectively. Hence, we omit details in order to avoid redundant arguments.



\subsection{Proof of Theorem \ref{thm:C1a-ini}}\label{sec:pf-thm-ini}

We can repeat the scaling argument used in the proof of Theorem \ref{thm:C1a}. The only difference is that we now also rescale the initial datum $g$. Since the modifications are now straightforward, we omit the details.

\subsection{Extension of  Propositions \ref{prop:C1a-relower}-\ref{prop:C1alower} to the full intrinsic cylinders}\label{extension} Combining the conclusions of Proposition \ref{prop:C1a-relower} and Proposition \ref{prop:C1alower}, valid in the lower intrinsic cylinders $(Q_r^-)$ only, with an application of Proposition \ref{prop:C1a-ini} with $g\equiv u$ for $\alpha\in (0,2/(p-2)]$ and $\alpha\in (0,\alpha_h)$, respectively, we can immediately state the following extensions of  Proposition \ref{prop:C1a-relower} and  Proposition \ref{prop:C1alower} to the full intrinsic cylinder.

\begin{proposition}\label{prop:C1a-re}
Under the assumptions in Proposition \ref{prop:C1a-relower},
there exists $\bar\lambda\equiv\bar\lambda(n,p,\nu,L,A,\alpha,\theta)$ such that
\begin{equation*}
|u - \ell | \leq \bar\lambda r\omega(r)\quad\text{in }Q_r(\bar\lambda\vp(r)),
\end{equation*}
for every $r\in(0,1)$.
\end{proposition}

\begin{proposition}\label{prop:C1a}
Under the assumption in Proposition \ref{prop:C1alower}, there exists $\bar\lambda\equiv\bar\lambda(n,p,\nu,L,A,\alpha)$ such that
\begin{equation}\label{eq:u-C1a}
|u - \ell | \leq \bar\lambda r\omega(r) \quad\text{in }Q_r(\bar\lambda\vp(r)),
\end{equation}
for every $r\in(0,1)$.
\end{proposition}

\subsection{Proof Theorem \ref{thm:C1a-re}}
The proof is analogous to the proof of Theorem \ref{thm:C1a} in Subsection \ref{pproffa}. In this case though we verify that Proposition \ref{prop:C1a-re} can be applied.


\section{Regularity across free boundaries}\label{finalproof}

This section is devoted to the proof of Theorem \ref{thm:C1a-int} and Theorem \ref{thm:C1a-gen}. These results will follow from Propositions \ref{prop:C1a-int}, \ref{prop:C1a} and \ref{prop:C1a-ini}. However, with the intrinsic geometry involved, it is not a trivial matter to put these results together. To describe the issue, let $(x_i,t_i)$ be two points, $i\in \{1,2\}$, with $t_1\neq t_2$, and let $u$ be a solution. Set
$$
d_i \equiv  \sup\left\{ r> 0: (x_j,t_j) \not\in (x_i,t_i) + Q_r(\vp_i(r)), i\neq j\right\},
$$
where
$$
\vp_i(r) \equiv \max\{\omega(r), |Du(x_i,t_i)| + s\}.
$$
Needless to say, $d_1$ and $d_2$ may not be equal, and these numbers may not even be comparable. Instead, the comparability of the intrinsic distance occurs only when comparability of the (non-)degeneracies, $|Du(x_1,t_1)| + s$ and $|Du(x_2,t_2)| + s$, can be ensured.

Keeping this in mind, we shall argue as follows. We first derive from Proposition \ref{prop:C1a} a uniform intrinsic $C^{1,\omega}$ approximation on the contact set. We then extend the approximation uniformly to each level surface that is equally distanced to the contact set with respect to the intrinsic geometry defined at the free boundary points. Upon extension, we show that the intrinsic distances are comparable.

We shall need the following elementary lemma.

\begin{lemma}\label{lem:alg}
Let $p > 2$, $s\geq 0$ and let $(\vec a,\vec b)$ be a pair of nonzero vectors. Then
$$
\int_0^1 ( |\vec a + \sigma \vec b|^2 + s^2 )^{\frac{p-3}{2}} \d\sigma \leq c\frac{(|\vec a| + s)^{p-2} + |\vec b|^{p-2}}{|\vec a| + s},
$$
for some constant $c \equiv c(p)$.
\end{lemma}

\begin{proof}
In this proof, $c$ is a positive constant depending at most on $p$, and which may vary upon each occurrence. Write $\mu\equiv |\vec a| + s$. Regardless of the value of $p$, if $|\vec b|\leq \frac{1}{2}\mu$, then
$$
\int_0^1 (|\vec a + \sigma \vec b|^2 + s^2 )^{\frac{p-3}{2}} \d\sigma \leq c \int_0^1 (\mu + \sigma |\vec b|)^{p-3}\d\sigma \leq c \mu^{p-3},
$$
so the assertion is proved in this case. Now let us consider the case $|\vec b| > \frac{1}{2}\mu$. In this case, if $p \geq 3$, then
\begin{align*}
\int_0^1  (|\vec a + \sigma \vec b|^2 + s^2 )^{\frac{p-3}{2}} \d\sigma \leq c \int_0^1 (\mu + \sigma |\vec b|)^{p-3}\d\sigma \leq \frac{c((\mu + |\vec b|)^{p-2} - \mu^{p-2})}{|\vec b|} \leq \frac{c(\mu^{p-2} + |\vec b|^{p-2})}{\mu},
\end{align*}
as desired. Furthermore, if $2 < p < 3$, then since $|\vec a + \sigma \vec b|^2 + s^2 \geq | |\vec a| - \sigma |\vec b| |^2 + s^2 \geq \frac{1}{2} | \mu - \sigma |\vec b||^2$,
\begin{align*}
\int_0^1  (|\vec a + \sigma \vec b|^2 + s^2 )^{\frac{p-3}{2}} \d\sigma \leq \int_0^1 |\mu - \sigma |\vec b||^{p-3}\d\sigma \leq  \frac{c(|\mu|^{p-2} + \vec |b|^{p-2})}{|\vec b|} \leq \frac{c(\mu^{p-2} + |\vec b|^{p-2})}{\mu},
\end{align*}
which verifies the assertion in this case.
\end{proof}


\subsection{Estimates across free boundaries}

We are going to present a detailed argument in the rescaled form.  Recall the H\"older exponent $\alpha_h$ from Theorem \ref{thm:grad -}.

\begin{proposition}\label{prop:C1a-gen}
Let $H$ be as in \eqref{generalH}, with $(a,s)$ as in \eqref{asp} and \eqref{eq:a-C2}, let $\omega$ verify \eqref{eq:omega prop} with $\alpha < \alpha_h$. Let $\psi$ be an obstacle such that $|D\psi| + s \leq 1$ in $Q_2$, and let $u$ be a solution to $\max\{Hu , \psi - u \} = 0$ in $Q_2$ such that $u(0,0) = 0$ and $|u| \leq 1$ in $Q_2$. Suppose that there exists, for each $(x_0,t_0)\in Q_1$,  a time-independent affine function $\tilde\ell$ such that
$$
 |\psi  - \tilde\ell| \leq r\omega(r)\quad\text{in }(x_0,t_0) + Q_r(\tilde\vp(r)),
$$
for every $r\in(0,1)$, with $\tilde\vp(r) \equiv \max\{\omega(r), |D\tilde\ell| + s\}$. Then there exists $\lambda\equiv\lambda(n,p,\nu,L,L',A,\alpha)$, such that for every $(x_0,t_0) \in Q_1$ there is a time-independent affine function $\ell$ such that
\begin{equation}\label{eq:C1a-gen}
| u - \ell  | \leq \lambda r\omega(r)\quad \text{in }(x_0,t_0) + Q_r(\lambda\vp(r)),
\end{equation}
for every $r\in(0,1)$, with $\vp(r)\equiv \max\{\omega(r), |D\ell| + s\}$.
\end{proposition}

\begin{proof}
Our starting point is the analysis at contact points. Let $(\tilde x,\tilde t)\in Q_1\cap \{ u = \psi\}$. Note that
$$
(\tilde x,\tilde t) + Q_1\subset Q_2.
$$
Therefore, Proposition \ref{prop:C1a} (after an obvious translation) yields \eqref{eq:C1a-gen} at $(\tilde x,\tilde t)$. In particular, the affine function $\ell$ and the (non-)degeneracy constant $\mu$ coincide with $\tilde\ell$ and $\tilde\mu$, respectively. More specifically,
\begin{equation}\label{eq:u-l'}
| u - \tilde\ell | \leq \lambda r\omega(r)\quad\text{in }(\tilde x,\tilde t) + Q_r(\lambda\tilde\vp(r)),
\end{equation}
for every $r\in(0,1)$. Note that $\tilde\ell$ is a time-independent affine function with $\tilde\ell(\tilde x) = \psi(\tilde x,\tilde t)$ and $D\tilde\ell = D\psi(\tilde x,\tilde t)$.

Now let $\varrho$ be given, with $0<\varrho<\frac{1}{6}$, and choose an arbitrary point
\begin{equation}\label{eq:xt}
(x_0,t_0) \in (\tilde x,\tilde t) + \partial Q_{4\varrho}(\lambda\tilde\vp(4\varrho)).
\end{equation}
We claim that for every $r\in(0,\frac{1}{6})$,
\begin{equation}\label{eq:u-l-claim}
 |u - \ell| \leq 3\bar\lambda r\omega(r)\quad\text{in }(x_0,t_0) + Q_r(\bar\lambda\vp(r)),
\end{equation}
for some constant $\bar\lambda\equiv \bar\lambda(n,p,\nu,L,L',A,\alpha)$, $\bar\lambda > 2\lambda$.

Let $\e$ be a sufficiently small constant to be determined later by $n$, $p$, $\nu$, $L$, $L'$, $A$ and $\alpha$ only. We divide our proof of the claim in \eqref{eq:u-l-claim} into two cases.

\begin{case}
$\omega(\varrho) \leq \e (|D\tilde\ell| + s)$.
\end{case}

In this case, we can linearize our problem. Write $\tilde\mu\equiv |D\tilde\ell| + s$. Since
\begin{equation}\label{eq:Q-Q'}
(x_0,t_0) + Q_{2\varrho}(\lambda\tilde\mu) \subset (\tilde x,\tilde t) + Q_{6\varrho}(\lambda\tilde\mu),
\end{equation}
it follows from \eqref{eq:u-l'} and that $\omega(\varrho)\leq \e\tilde\mu < \tilde\mu$, that
\begin{equation}\label{eq:Du-sup}
|Du| \leq c_0 \lambda \tilde\mu \quad\text{in }(x_0,t_0) + Q_\varrho(\lambda\tilde\mu),
\end{equation}
for some $c_0\equiv c_0(n,p,\nu,L)$, $c_0\geq 1$, and that
\begin{equation}\label{eq:Du-osc}
|Du - Du(x_0,t_0)| \leq c_0\lambda \tilde\mu \sigma^{\alpha_h}\quad\text{in }(x_0,t_0) + Q_{\sigma\varrho}(\lambda\tilde\mu),
\end{equation}
for every $\sigma \in (0,1)$, with $\alpha_h$ as in Theorem \ref{thm:grad -}. We introduce the change of variables
$$
(y,\tau) \equiv ( {\varrho^{-1}}{(x - x_0)}, {{\varrho^{-2}}(\lambda\tilde\mu)^{p-2}(t - t_0)})\quad\text{and}\quad w (y,\tau) \equiv \frac{(u-\tilde\ell)(x,t)}{\lambda\varrho\tilde\mu}.
$$
Using that $Hu = H\tilde\ell = 0$ in $(x_0,t_0) + Q_{2\varrho}(\lambda\tilde\mu)$, we see that
\begin{equation}\label{eq:u-l'-pde}
\ddiv (b Dw ) = w_\tau \quad\text{in }Q_1,
\end{equation}
where $b$ is the $(n\times n)$-dimensional matrix-valued coefficient,
$$
b(y,\tau)\equiv \frac{1}{(\lambda\tilde\mu)^{p-2}}\int_0^1 \partial_z a (\eta Du(x,t) + (1-\eta)D\tilde\ell)\,\d\eta.
$$
Since $\tilde\mu \equiv |D\tilde\ell| + s$, we can deduce from \eqref{eq:Du-sup}, as we did in \eqref{eq:b-ellip-up} and \eqref{eq:b-ellip-low}, that
\begin{equation}\label{eq:b-ellip}
\frac{c_1\nu}{\lambda^{p-2}}|\xi|^2 \leq \langle b\xi,\xi\rangle \leq c_2 c_0^{p-2} L|\xi|^2\quad\text{in }Q_1,\quad\forall  \xi\in\R^n\setminus\{0\},
\end{equation}
for some constants $c_i \equiv c_i(n,p)$, $i\in\{1,2\}$, with $0<c_1 < 1< c_2$. Moreover, by the assumption in \eqref{eq:a-C2} on $a$, and \eqref{eq:Du-osc}, we can derive by linearizing of the coefficient $b$, that
\begin{equation}\label{eq:b-osc}
| b - b(0,0) | \leq  \frac{c_0\lambda L'B}{(\lambda\tilde\mu)^{p-3}}\sigma^{\alpha_h}\quad\text{in }Q_\sigma,
\end{equation}
for every $\sigma\in(0,1)$, where we have set
$$
B \equiv \sup_{(x_0,t_0) + Q_\varrho(\lambda\tilde\mu)} \int_0^1 \int_0^1  (|(1-\eta)D\tilde\ell +\eta D ( \eta' u+ (1-\eta')u(x_0,t_0))|^2 + s^2)^{\frac{p-3}{2}}\,\d\eta\,\d\eta'.
$$
By Lemma \ref{lem:alg} and \eqref{eq:Du-sup} (and recalling that $\tilde\mu \equiv |D\tilde\ell| + s$),
\begin{equation}\label{eq:B}
B \leq \frac{c_4}{\tilde\mu} \left[ \tilde\mu^{p-2} + \sup_{(x_0,t_0) + Q_\varrho(\lambda\tilde\mu)} |Du|^{p-2} \right] \leq 2c_4 (c_0\lambda)^{p-2}\tilde\mu^{p-3},
\end{equation}
for some $c_4\equiv c_4(p)$. Hence, combining \eqref{eq:b-osc} and \eqref{eq:B}, we obtain that
\begin{equation}\label{eq:b-osc-re}
|b - b(0,0)| \leq  2c_4c_0^{p-1}\lambda L' \sigma^{\alpha_h}\quad\text{in }Q_\sigma,
\end{equation}
whenever $\sigma\in(0,1)$. In view of \eqref{eq:b-ellip} and \eqref{eq:b-osc-re}, we can invoke the pointwise $C^{1,\alpha_h}$ approximation for uniformly parabolic linear problems, see \cite{Lieberman_book}. Since \eqref{eq:u-l'}, along with \eqref{eq:Q-Q'} and \eqref{eq:omega prop}, gives us
\begin{equation}\label{eq:w-l''-re}
 |w| \leq \frac{6^{1+\alpha} A\omega(\varrho)}{\tilde\mu}\quad\text{in }Q_1,
\end{equation}
the parabolic regularity theory, see \cite{Lieberman_book},  yields a time-independent affine function $\tilde\ell'$ such that
\begin{equation}\label{eq:w-l''}
\begin{aligned}
 | w - \tilde\ell'| &\leq \frac{6^{1+\alpha} A\omega(\varrho)}{\tilde\mu}c_5\sigma^{1 + \alpha_h}\quad\text{in }Q_\sigma,
\end{aligned}
\end{equation}
for every $\sigma\in(0,1)$, and for some constant $c_5\equiv c_5(n,p,\nu,L,L',\lambda) \equiv c_5(n,p,\nu,L,L')$, $c_5 \geq 1$.

At this point, we select
$$
\ell(x) \equiv \tilde\ell(x) + \lambda\varrho\tilde\mu \tilde\ell'\left(\frac{x - x_0}{\varrho}\right),
$$
as the approximating affine function $\ell$ for $u$ at $(x_0,t_0)$. Then \eqref{eq:w-l''} can be rephrased in the original scale as
\begin{equation}\label{eq:u-l-claim-re}
 | u - \ell | \leq 6^{1+\alpha} A \lambda c_5 \sigma^{1+\alpha_h} \varrho\omega(\varrho) \leq 6^{1+\alpha} A^2 \lambda c_5\sigma \varrho\omega(\sigma\varrho)\quad\text{in }(x_0,t_0) + Q_{\sigma\varrho}(\lambda\tilde\mu),
\end{equation}
for all $\sigma\in(0,1)$, where the rightmost inequality is deduced from \eqref{eq:omega prop} and \eqref{eq:alpha}. By \eqref{eq:w-l''-re} and \eqref{eq:w-l''}, we deduce that
\begin{equation}\label{eq:Dl-Dl'}
|D\ell - D\tilde\ell| \leq 2\cdot 6^{1+\alpha} A\lambda c_5 \omega(\varrho).
\end{equation}
Now writing $\mu\equiv |D\ell| + s$, and recalling that we in this case assume $\omega(\varrho) \leq \e\tilde\mu$, it follows from \eqref{eq:Dl-Dl'} that
\begin{equation}\label{eq:m-m'}
\left| \frac{\mu}{\tilde\mu} - 1 \right|  \leq \frac{1}{2},
\end{equation}
provided that $4\cdot 6^{1+\alpha} \lambda c_5\e \leq 1$. The smallness condition for $\e$ is determined solely by $A$, $\alpha$, $\lambda$ and $c_5$, in fact,  tracking the dependence of these constants, we can choose $\e\equiv \e(n,p,\nu,L,L',A,\alpha)$, $0< \e < \frac{1}{144}$. With \eqref{eq:m-m'} at hand, we have
$$
Q_{\sigma\varrho}(2\lambda\mu) \subset Q_{\sigma\varrho}(\lambda\tilde\mu),
$$
which in \eqref{eq:u-l-claim-re} yields that
\begin{equation}\label{eq:u-l-claim-re2}
 |u - \ell| \leq \bar\lambda r \omega(r)\quad\text{in }(x_0,t_0) + Q_r(\lambda\mu),
\end{equation}
for all $r\in(0,\varrho)$, for some constant $\bar\lambda \equiv \bar\lambda(n,p,\nu,L,L',A,\alpha)$, with $\bar\lambda > 2\cdot 6^{1+\alpha} A^2\lambda c_5$. With such a choice of $\bar\lambda$, we can use \eqref{eq:Dl-Dl'} and \eqref{eq:Q-Q'} to deduce from \eqref{eq:u-l'} that
\begin{equation}\label{eq:u-l-claim-re3}
| u - \ell| \leq \sup_{(\tilde x,\tilde t) + Q_{6r}(\lambda\tilde\vp(r))} | u - \tilde\ell | + \sup_{B_{6r}(x_0)} |\ell - \tilde\ell| \leq 3\bar\lambda r \omega(r)\quad\text{in }(x_0,t_0) + Q_r(\bar\lambda\vp(r)),
\end{equation}
for all $r\in [\varrho,\frac{1}{6})$, with $\vp(r)\equiv \max\{\omega(r), |D\ell| + s \}$. Combining \eqref{eq:u-l-claim-re2} with \eqref{eq:u-l-claim-re3}, we verify our claim in \eqref{eq:u-l-claim} for Case 1.

\begin{case}
$\omega(\varrho) \geq \e(|D\tilde\ell| + s)$.
\end{case}

As for this case, we use Proposition \ref{prop:C1a-int}. By \eqref{eq:u-l'}, we obtain
$$
|u - u(x_0,t_0)| \leq \frac{\lambda}{\e}\varrho\omega(\varrho)\quad\text{in }(x_0,t_0) + Q_\varrho(\lambda\e^{-1}\omega(\varrho)).
$$
Hence, Proposition \ref{prop:C1a-int} along with suitable rescaling and \eqref{eq:omega prop} yields
\begin{equation}\label{eq:u-l-re3}
\begin{aligned}
 |u - \ell|  \leq\bar\lambda \sigma\varrho\omega(\sigma\varrho)\quad\text{in }(x_0,t_0) + Q_{\sigma\varrho}(\bar\lambda\vp(r))), \quad\forall  \sigma\in(0,1],
\end{aligned}
\end{equation}
with $\ell$ being the time-independent affine function such that $\ell(x_0) = u(x_0,t_0)$ and $D\ell = Du(x_0,t_0)$, $\vp(r)\equiv \max\{\omega(r), |D\ell| + s \}$, $\bar\lambda \equiv c_6\lambda A/\e$ and $c_6 \equiv c_6(n,p,\nu,L,\alpha)$, $c_2\geq 1$. In comparison with \eqref{eq:u-l'} (at $r = \varrho$), we deduce as in \eqref{eq:Dl-Dl'} and with $A\geq 1$, that
\begin{equation}\label{eq:Dl-Dl'-re}
|D\ell - D\tilde\ell| \leq 2\bar\lambda\omega(\varrho).
\end{equation}
On the other hand, $\omega(\varrho)\geq \e(|D\tilde\ell| + s)$ implies that $\vp(r) \geq \e \tilde\vp(r)$, whenever $r\geq\varrho$, where $\tilde\vp(r) \equiv \max\{\omega(r), |D\tilde\ell| + s\}$. Since $\bar\lambda \geq\lambda/\e$,
\begin{equation}\label{eq:Q-Q'-re}
(x_0,t_0) + Q_r(\bar\lambda\vp(r)) \subset (\tilde x,\tilde t) + Q_{6r}(\lambda\tilde\vp(r)),
\end{equation}
whenever $r \geq \varrho$. With \eqref{eq:u-l'} and \eqref{eq:Dl-Dl'} at hand, we use \eqref{eq:Q-Q'-re} to derive, similarly as in the deduction of \eqref{eq:u-l-claim-re3}, that
\begin{equation}\label{eq:u-l-re4}
|u - \ell| \leq 3\bar\lambda r\omega(r)\quad\text{in }(x_0,t_0) + Q_r(\bar\lambda\vp(r)),
\end{equation}
for every $r\in(\varrho,\frac{1}{6})$. Hence, Case 2 is also settled.

Finally, we are left with the task of checking the  assertion in \eqref{eq:u-l-claim} but at any point $(x_0,t_0)$ that does not satisfy \eqref{eq:xt} for any contact point $(\tilde x,\tilde t)$ and for any $\varrho\in(0,\frac{1}{6})$. If such a point $(x_0,t_0)$ exists, then since $\tilde\vp(\frac{2}{3}) \equiv\tilde\vp(\frac{2}{3};\tilde x,\tilde t)\leq 1$ (because we assume in the statement of this proposition that $|D\psi(\tilde x,\tilde t)| + s \leq 1$ and $\omega(\frac{2}{3}) \leq \omega(1) = 1$) and $(x_0,t_0)\in Q_1$, we must have $u > \psi$ in $(x_0,t_0) + Q_{2/3}(\bar\mu) \subset Q_2$ for some large $\bar\mu\equiv\bar\mu(n,p,\nu,L,L',A,\alpha)$. Thus,
$$
H u = 0 \quad\text{in }(x_0,t_0) + Q_{2/3}(\bar\mu),
$$
and the desired estimate follows now from Proposition \ref{prop:C1a-int}. We omit the details for this last part, as it repeats some of our argument above.
\end{proof}

\begin{remark}\label{rem:C1a-gen}
The assumption in \eqref{eq:a-C2} is used to treat Case 1 above, which concerns the situation where the point on the free boundary is less degenerate. However, following the algebraic manipulations in \cite[Page 724]{Lin}, one can obtain, in place of \eqref{eq:C1a-gen}, that
$$
|u - \ell| \leq \lambda r\omega(r^\beta),
$$
for some $\beta\equiv \beta(n,p,\nu,L)$, $\beta \in (0,1)$. We leave the details to the interested reader.
\end{remark}

By a similar argument, we obtain a rescaled version of our result for the regularity at the initial layer/state. Since the proof is essentially the same, and there is no ambiguity in gluing the estimate in the forward cylinders (Proposition \ref{prop:C1a-ini}) with the full cylinders (Proposition \ref{prop:C1a-gen}), we shall omit the details and we simply state the result as follows.

\begin{proposition}\label{prop:C1a-gen-ini}
Let $H$ be as in \eqref{generalH}, with $(a,s)$ as in \eqref{asp} and \eqref{eq:a-C2},  let $\omega$ verify \eqref{eq:omega prop} with $\alpha\in (0,\alpha_h)$. Let $\psi$ be an obstacle such that $|D\psi| + s \leq 1$ in $Q_2^+\cup\partial_{\rm P}Q_2^+$, let $g$ be an initial datum such that $|Dg| + s \leq 1$ and $g\geq \psi(\cdot,0)$ in $B_2$. Let $u$ be a solution to
$$
\begin{dcases}
\max\{ Hu , \psi- u\} = 0 &\text{in }Q_2^+,\\
u = g &\text{on } B_2\times \{0\},
\end{dcases}
$$
such that $|u| \leq 1$ in $Q_2^+\cup\partial_{\rm P}Q_2^+$ and $u(0,0) = g(0) = 0$. Assume that there are, for each $(x_0,t_0)\in Q_1^+\cup(B_1\times\{0\})$, time-independent affine functions $\tilde \ell'$ and $\tilde\ell$ such that
$$
|g - \tilde \ell'| \leq r \omega(r)\quad\text{in }B_r(x_0),\quad |\psi  - \tilde\ell| \leq r\omega(r)\quad\text{in }(x_0,t_0) + Q_r^\pm(\tilde\vp(r)),
$$
whenever $(x_0,t_0) + Q_r^\pm(\tilde\vp(r))\subset Q_1^+$, with $\tilde\vp(r)\equiv \max\{\omega(r), |D\tilde\ell| +s \}$. Then there exists $\lambda\equiv\lambda(n,p,\nu,L,L',A,\alpha)$ such that for every $(x_0,t_0) \in Q_1^+$, there exists a time-independent affine function $\ell$ such that
\begin{equation}\label{eq:C1a-gen-ini}
| u - \ell  | \leq \lambda r\omega(r)\quad \text{in }(x_0,t_0) + Q_r^\pm(\lambda\vp(r)),
\end{equation}
whenever $(x_0,t_0) + Q_r^\pm(\tilde\vp(r))\subset Q_1^+$, with $\vp(r)\equiv\max\{\omega(r), |D\ell| + s\}$.
\end{proposition}


\subsection{Proofs of Theorems \ref{thm:C1a-int} and \ref{thm:C1a-gen}}

With Proposition \ref{prop:C1a-gen} at hand, Theorem \ref{thm:C1a-int} can be proved via a similar rescaling argument as the one used and detailed in the proofs of Theorems \ref{thm:C1a} and \ref{thm:C1a-re}. As for Theorem \ref{thm:C1a-gen}, we can simply apply Propositions \ref{prop:C1a-gen} and \ref{prop:C1a-gen-ini} at every points in the given subdomain $\O'$ of $\O$. Hence, the proofs are complete.\\

\noindent
{\bf Acknowledgement.}  Originally this was supposed to be a joint paper with Tuomo Kuusi and Giuseppe Mingione. In fact, after the completion of \cite{KMNobs}, Tuomo Kuusi in 2015 wrote a first draft for this continuation. Though our paper is much different compared to that draft, we learned a
lot from Tuomo Kuusi's draft.  We would have preferred to make this a joint paper. We acknowledge the generosity of Tuomo Kuusi and Giuseppe Mingione for letting us go ahead and finish this paper.

\end{document}